\documentclass[12pt]{amsart}

\usepackage[utf8]{inputenc}
\usepackage[T1]{fontenc}
\usepackage{amsmath,amssymb,amsthm}
\usepackage[shortlabels]{enumitem}
\usepackage{float}
\usepackage{subfiles}
\usepackage{hyperref}
\usepackage{bbm}
\usepackage{xcolor}

\usepackage{apptools}
\usepackage[titletoc]{appendix}

\theoremstyle{plain}
\newtheorem{theorem}{Theorem}[section]
\newtheorem*{theorem*}{Theorem}
\newtheorem{lemma}[theorem]{Lemma}
\newtheorem{proposition}[theorem]{Proposition}

\newtheorem*{properties*}{Properties}
\newtheorem*{lemma*}{Lemma}

\theoremstyle{definition}
\newtheorem*{definition*}{Definition}
\newtheorem{definition}{Definition}
\newtheorem*{example*}{Example}

\newtheorem*{remark*}{Remark}

\newcommand{\R}{\mathbb{R}}
\newcommand{\N}{\mathbb{N}}

\title[Construction and properties for the Green's function]{Construction and properties for the Green's function with Neumann boundary condition}
\author{Antoine Bricmont}
\date{}
\begin{document}
\begin{abstract}
This article addresses the construction and analysis of the Green's function for the Neumann boundary value problem associated with the operator \(-\Delta + a\) on a smooth bounded domain \(\Omega \subset \mathbb{R}^N\) (\(N \geq 3\)) with $a\in L^\infty(\Omega)$. Under the assumption that \(-\Delta + a\) is coercive, we obtain the existence, uniqueness, and qualitative properties of the Green's function \(G(x,y)\). The Green's function \(G(x,y)\) is constructed explicitly, satisfying pointwise estimates and derivative estimates near the singularity. Also, near the boundary of \(\Omega\), \(G\) is compared to the Green's function of the laplacian, with pointwise estimates. Other properties, like symmetry and positivity among other things, are established.
\end{abstract}
\maketitle
\section{Introduction}

We consider the following problem
\begin{align}
	\begin{cases}\label{problem}
		-\Delta  u +  a  u &=  f  \text{ in } \Omega, \\
		\partial_\nu  u &= 0   \text{ on } \partial \Omega,
	\end{cases}
\end{align}
where $\Omega$ is a smooth bounded open subset of $\R^N$ with $N\geq 3$, $\Delta:=\sum\limits_{i=1}^N \partial^2_{ii} $, $u \in C^2(\overline\Omega)$, $a\in L^\infty(\Omega)$, $f\in C^{0,\alpha}(\overline\Omega)$, and $\nu(x)$ is the unit outward vector at $x\in \partial \Omega$. We also assume that $-\Delta + a$ is coercive, i.e.\ for all $u\in H^1(\Omega)$, there exists $C>0$ such that
$$
\int_{\Omega}|\nabla u |^2 \,+ a u^2 \,dx \geq C \int_{\Omega}|\nabla u |^2 \,+ u^2 \,dx.
$$
The purpose of this article is to construct and give some properties of the Green's function associated to (\ref{problem}). Green's functions provide explicit solution formulas for partial differential equations but also encode information about differential operators, domain geometry, and boundary conditions. 

Let $S\subseteq \R^N$, we recall that $Diag(S):= \{ (x,x) | x \in S \}$, and give the definition of a Green's function for problem (\ref{problem}).
\begin{definition}[Green's function]
	Let $a\in L^\infty(\Omega)$ and $G: \overline{\Omega} \times \overline\Omega \setminus Diag( \overline\Omega)\rightarrow \R$. We say that $G$ is a Green's function of $-\Delta + a$ with Neumann condition on the boundary of $\Omega$ if, denoting $G_x:= G(x,.)$, the following affirmations are satisfied for all $x \in \Omega$:
	\begin{enumerate}[(i)]
		\item $G_x \in L^1(\Omega )$; 
		\item for all $\phi \in C^2(\overline\Omega)$, we have $$ \int_\Omega G_x(-\Delta \phi + a \phi ) \, dy = \phi (x) - \int_{\partial \Omega} G_x \partial_\nu \phi \,d\sigma,$$
	where $d\sigma$ is the surface element of $\partial\Omega$.	
	\end{enumerate}
\end{definition} 

In this paper, we construct a unique Green's function $G$ for problem (\ref{problem}), give several pointwise estimates for $G$ and its derivatives, and give some key properties that $G$ possesses.
Our main result is the following theorem.
\begin{theorem}\label{theorem_main}
	Let $a\in L^\infty(\Omega)$. \begin{enumerate}[(i)]
		\item There exists a unique Green's function $G$ associated with problem (\ref{problem}).
		\item There exists $C>0$, such that 
	\begin{align*} 
	| G(x,y) |+|x-y|| \nabla_y G(x,y) |\leq \frac{C}{|x-y|^{N-2}}
	\end{align*}
	for all $x,y \in  \overline\Omega$, $x\neq y$.
	\item Let $\Gamma$ be the Green's function of operator $-\Delta$ on $\R^N$. Let also $x\in K\subset \Omega$ a compact set. Then there exists $C=C(K) >0$, such that 
	\begin{align*}
	\left|G(x,y)-\Gamma(x,y)\right|&\leq C|x-y|^{3-N},\\
	\left| \nabla_y G(x,y)-\nabla \Gamma(x,y)\right|&\leq C \begin{cases}
		|\ln|x-y|| &\text{ if }N=3\\
		|x-y|^{3-N} &\text{ if }N>3
	\end{cases},
\end{align*}
for all $y\in \overline\Omega$, $y\neq x$.
	\item Let $\Gamma_\nu$ be the Green's function of operator $-\Delta$ with Neumann boundary condition on $\R_+^N$. Let $x_0\in \partial \Omega$ and $\varphi$ be a chart as in Definition \ref{smooth_subset} such that $\varphi(0)=x_0$, then there exists $C>0$ and $R>0$ such that 
	 \begin{align*}
	 	\left|G(\varphi(x),\varphi(y))-\Gamma_\nu(x,y)\right|\leq C\begin{cases}
		|\ln|x-y|| & \text{ if }N=3 \\
		|x-y|^{3-N} & \text{ if } N>3,
	\end{cases}\\
	\left| \partial_{y_i} G(\varphi(x),\varphi(y))-\partial_{y_i} \Gamma_\nu(x,y)\right|\leq C \begin{cases}
		\dfrac{|\ln|x-y||}{|x-y|} & \text{ if }N=3 \\
		|x-y|^{2-N} & \text{ if } N>3,
		\end{cases}
	 \end{align*}
for all $x,y \in B(0,R)\cap\R^N_+$, $y\neq x$.
	\end{enumerate}
\end{theorem}
The functions $\Gamma$ and $\Gamma_\nu$ are explicit and given in section \ref{section_construction} and \ref{section_partial_Omega} respectively.

We prove Theorem~\ref{theorem_main} by adapting the approach in \cite{Druet-Robert-Wei} 	and \cite{Robert-Dirichlet} (see also \cite{Aubin}, and \cite{Aubin-book-1998}) and by explicitly constructing $G$ via an iteration procedure. The main novelty of our work is the point (iv) of Theorem~\ref{theorem_main} where we obtain quantitative estimates on $G$ near the boundary. All the results that we prove in this work  will be used in a future article. For other (non-exhaustive) examples of constructions and properties of Green's functions for related problems, we refer to \cite{Druet-Premoselli-2014}, \cite{Fine-Premoselli-2020}, \cite{Ghoussoub-Mazumdar-Robert-2023}, and  \cite{Premoselli-2016}.

\section{Construction of the Green’s function}
\label{section_construction}
All this section follows closely \cite{Druet-Robert-Wei} and \cite{Robert-Dirichlet}. We first recall a very classical definition.
\begin{definition}\label{smooth_subset}
	We say that $\Omega$ is smooth if for all $x \in \partial \Omega$, there exists $R(x)=R > 0$, $U_x$ an open neighborhood of $x$ in $\R^N$, and $\hat\varphi : B(0,R) \rightarrow U_x$, where $B(0,R)$ is the ball centered in $0$ and 
of radius $R$ in $\R^N$, such that
\begin{enumerate}[(i)]
		\item $\hat\varphi$ is a $C^\infty$-diffeomorphism;
		\item $\hat\varphi (0)=x$;
		\item $\hat\varphi (B(0,R)^+)=\hat\varphi (B(0,R)) \cap \Omega $, where $B(0,R)^+:=B(0,R)\cap \{x_N>0\}$;
		\item $\hat\varphi (B(0,R) \cap \{x_N=0\} )=\hat\varphi (B(0,R)) \cap \partial\Omega$.
	\end{enumerate}
\end{definition}

We set
$$
\Gamma(x,y):= \frac{1}{(N-2)\omega_{N-1}|x-y|^{N-2}}
$$
for all $(x,y) \in \R^N \times \R^N \setminus Diag(\R^N)$ where $\omega_{N-1}$ is the volume of the unit sphere of $\R^N$. We recall certain properties that $\Gamma$ satisfies:\\ 
	$\Gamma \in C^\infty\left( \R^N \times \R^N \setminus Diag(\R^N) \right) $ and is such that
	\begin{enumerate}[(i)]
		\item \label{prop_Heq1} for all $x,y \in \Omega$ such that $x\neq y$, $$ |x-y|^{N-2} |\Gamma(x,y)| + |x-y|^{N-1} |\nabla \Gamma(x,y)| \leq \frac{N-1}{(N-2)\omega_{N-1}} ; $$
		\item  for all $x \in \Omega$, $$ -\Delta \Gamma_x = \delta_x \text{ in }D'(\Omega);$$
		\item \label{prop_Heq3} for all $\phi \in C^2(\R^N)$ and for all $x \in \R^N$, \begin{align*}
				 \int_\Omega \Gamma_x (-\Delta \phi) \, dy = \phi (x) + \int_{\partial \Omega} (\phi \partial_\nu \Gamma_x - \partial_\nu \phi \Gamma_x) \,d\sigma. 
			\end{align*} 
	\end{enumerate}
The construction of the Green's function relies on two propositions. The first one is the following.
\begin{proposition}\label{construct}
	Let $a \in L^\infty(\Omega)$. Then, there exists a function $\hat G \in L^\infty_{loc}\left(\overline \Omega \times \overline \Omega \setminus Diag(\overline \Omega) \right) $ such that for all $x\in  \overline\Omega $, we have $\hat G_x \in C^{1,\alpha}_{loc}\left( \overline \Omega \setminus\{x \} \right)$ for all $\alpha \in (0,1) $ and for all $\phi \in C^2(\overline  \Omega )$, we have 
	\begin{align*}
		\int_\Omega \hat G_x(-\Delta \phi + a \phi ) \, dy = \phi (x) + \int_{\partial  \Omega} (\phi \partial_\nu \hat G_x - \hat G_x \partial_\nu \phi) \,d\sigma.
	\end{align*}
	Furthermore, there exists $C>0$, such that $ |\hat G(x,y) |\leq \frac{C}{|x-y|^{N-2}}$ for all $x,y \in  \overline\Omega$, $x\neq y$.
\end{proposition}
\begin{proof}
	We define
	\begin{align*}
		&\Gamma_1 (x,y) := -a(y) \Gamma(x,y),
	\end{align*}
	and for all integers $i\geq 1$
	\begin{align*}
		& \Gamma_{i+1} (x,y) := \int_\Omega \Gamma_i(x,z) \Gamma_1 (z,y) \,dz , 
	\end{align*}
	for all $x,y \in \overline \Omega, x\neq y$. So that, using Giraud's lemma (Proposition~\ref{giraud}), $\Gamma_i \in L_{loc}^\infty\left(\overline \Omega \times \overline \Omega \setminus Diag(\overline \Omega) \right)$, and there exists $C(i)>0$ such that
	$$ |\Gamma_i(x,y)|\leq C(i) \begin{cases}
		|x-y|^{2i-N} &\text{ if }i<\frac{N}{2} \\
		1+\log|x-y| &\text{ if }i=\frac{N}{2} \\
		1 &\text{ if }i>\frac{N}{2}
	\end{cases}, $$
	 for all $i\geq 1$ and all $x,y \in \overline \Omega, x\neq y$. Furthermore, when $i>\frac{N}{2}$, $\Gamma_i \in C^0(\overline\Omega \times \overline\Omega)$.
	
	Now we take $k$ the largest integer such that $k\leq\frac{N}{2}$, and we define
	$$
	G'(x,y) := \Gamma(x,y) + \sum_{i=1}^{k}\int_\Omega \Gamma_i(x,z)\Gamma(z,y)\,dz
	$$
	for all $x,y \in \overline \Omega, x\neq y$.
	Thanks to Giraud's lemma (Proposition~\ref{giraud}), we have 
	\begin{align}
	|G'(x,y)|\leq C|x-y|^{2-N}
	\end{align}
	for all $x,y \in \overline \Omega, x\neq y$, $G' \in C^0(\overline\Omega \times \overline\Omega\setminus Diag(\overline \Omega))$, and there exists $C>0$ such that
	\begin{align}
		|G'(x,y)-\Gamma(x,y)| \leq C|x-y|^{3-N}
	\end{align}
	for all $x,y \in \overline \Omega, x\neq y$.
	
	We now compute $-\Delta G'_x + a G'_x$ in the sense of distributions. Thanks to the estimates on $\Gamma_i$ above, we can apply Fubini's theorem. Let $\phi \in C^2(\overline  \Omega)$, we have 
	\begin{align*}
		&\int_\Omega (-\Delta \phi + a \phi) G'_x \,dy =\int_\Omega (-\Delta \phi+a \phi) \Gamma_x \,dy \\
		+& \sum_{i=1}^k \int_\Omega\int_{\Omega } (-\Delta \phi(y) + a(y) \phi(y)) \Gamma_i(x,z)\Gamma(z,y)\,dy \, dz.
	\end{align*} 
	Using property \ref{prop_Heq3} of $\Gamma$ and Fubini again, we get
	\begin{align*}
		&\int_\Omega (-\Delta \phi + a \phi) G'_x \,dy \\
		= \, &\phi (x) + \int_{\partial \Omega} (\phi \partial_\nu \Gamma_x - \partial_\nu \phi \Gamma_x) \,d\sigma + \int_\Omega  a \, \Gamma_x  \phi \, dy \\
		&+ \sum_{i=1}^k \int_\Omega\int_{\Omega } a(y) \phi(y)\, \Gamma_i(x,z)\Gamma(z,y)\,dz \, dy \\ 
		&+  \sum_{i=1}^k\int_\Omega \Gamma_i(x,z) \left[  \phi (z) + \int_{\partial \Omega} (\phi \partial_\nu \Gamma_z- \partial_\nu \phi \Gamma_z) \,d\sigma \right] dz \\
		= \,  & \phi (x) + \int_\Omega a \Gamma_x \phi \, dy  +  \sum_{i=1}^k\int_\Omega\int_{\Omega } a(y) \phi(y)\, \Gamma_i(x,z)\Gamma(z,y)\,dz \, dy \\ 
		&+\sum_{i=1}^k\int_\Omega \Gamma_i(x,z) \phi(z)\,dz  \, + \int_{\partial \Omega} (\phi \partial_\nu G'_x - \partial_\nu \phi G'_x) \,d\sigma, 
		\end{align*} 
		so that by definition of the $\Gamma_i$'s, we obtain
		\begin{align*}
			&\int_\Omega (-\Delta \phi + a \phi) G'_x \,dy \\
		= \, & \phi(x)+ \int_{\partial \Omega} (\phi \partial_\nu G'_x - \partial_\nu \phi G'_x) \,d\sigma - \int_\Omega\int_\Omega \Gamma_k(x,z)\Gamma_1(z,y)\phi(y) \,dz\,dy\\
		= \, & \phi(x)+ \int_{\partial \Omega} (\phi \partial_\nu G'_x - \partial_\nu \phi G'_x) \,d\sigma -\int_\Omega\Gamma_{k+1}(x,y)\phi(y)\,dy.
		\end{align*}	
		Thus, $(-\Delta + a) G_x'=\delta_x-(\Gamma_{k+1})_x$ in the sense of distributions. 
		By standard elliptic regularity \cite[Theorem 9.11]{Gilbarg-Trudinger} we can deduce that $G_x'\in C^{1,\alpha}_{loc}(\overline\Omega \setminus\{x\})$ for all $x \in \overline\Omega$ and all $\alpha \in (0,1)$.
	
	Now, by Theorem~\ref{thm_exist} there exists  $v_x \in H^1(\Omega )$, a unique weak solution of 
	$$
	\begin{cases}
		-\Delta v_x+a v_x &= (\Gamma_{k+1})_x  \text{ in } \Omega \\
		\partial_\nu v_x&=0   \text{ on }  \partial \Omega.
	\end{cases}
	$$
	so that by standard elliptic regularity \cite[Theorem 2.4.2.7 and Theorem 2.3.3.6]{Grisvard}, there exists $C>0$ independent of $x$ such that $\Vert v_x \Vert_{C^{1,\alpha}(\overline \Omega)}\leq C$ for all $\alpha \in (0,1)$, because $(\Gamma_{k+1})_x$ is continuous on $\overline\Omega$ for all $x \in \overline \Omega$.
	
	We set $\hat G_x(y) = G_x'(y) + v_x(y)$, for all $x,y \in \overline \Omega, x\neq y$ and $\hat G(x,y)$ verifies
	$$
	\begin{cases}
		-\Delta \hat G_x+a \hat G_x &= 0  \text{ in } \Omega \setminus\{x\} \\
		\partial_\nu \hat G_x &=\partial_\nu G_x'  \text{ on }  \partial \Omega.
	\end{cases}
	$$
\end{proof}

In order to complete the construction of the Green's function and to prove estimates up to the boundary, we introduce a classical way to extent a function outside of $\Omega$.

Let $x_0=(x_0',(x_0)_N)\in \partial \Omega$, and $R>0$, we set 
\begin{align*}
\rho(x'):=\hat\varphi_N(x',0)
\end{align*}
where $\hat \varphi$ is a chart as in Definition \ref{smooth_subset}, so that $\rho : \{x'\in \R^{N-1} \mid |x'|<R\}\rightarrow \R$ is a smooth function such that 
$$
\begin{cases}
	\Omega \cap B(x_0,R) = \{(x',x_N) \in B(x_0,R) \mid x_N>\rho(x') \},\\
	\partial\Omega \cap B(x_0,R) = \{(x',x_N) \in B(x_0,R) \mid x_N=\rho(x') \}.	
\end{cases}
$$

Up to an isometry, we can choose the orthogonal coordinate system in $\R^N$ in such a way that the first $N-1$ coordinate directions lie in the tangent space to $\partial \Omega$ at $x_0$ and the $N$th coordinate direction lie on the inward normal direction at $x_0$. Thus, we can take $\rho(x'_0)=0, \nabla \rho(x'_0)=0$.

We observe that an outward normal vector is given by 
\begin{align}\label{vec_normal}
n =(\nabla \rho, -1).
\end{align} 

\begin{definition}\label{varphi_pi}
Let $V\subset \R^N$ be an open set containing the origin, we define 
\begin{align*}
&\varphi_{x_0}:V\rightarrow \R^N:(x',x_N) \rightarrow (x',\rho(x'))-x_N  n(x',\rho(x')),\\
\text{and }&\\
&\pi:\R^N\rightarrow \R^N:(x',x_N) \rightarrow (x',|x_N|). 
\end{align*}
\end{definition}
We notice that $\varphi_{x_0}(x_0)=x_0$, since $(x_0)_N=\rho(x_0)=0.$ In what follows, we might omit the subscript $x_0$ when it is not necessary in context.

Since the Jacobian of $\varphi$ is given by 
\begin{align*}
J_{\varphi}(x)&=\\
&\begin{pmatrix}
1-x_N\partial_1   n_1 & -x_N\partial_2  n_1 & \cdots & -x_N\partial_{N-1}n_1 &-  n_1 \\
-x_N\partial_1  n_2 & 1-x_N\partial_2   n_2 & \cdots & - x_N \partial_{N-1}n_2 &-  n_2 \\
\vdots & \vdots & & \vdots & \vdots  \\
-x_N\partial_1 n_{N-1} & -x_N\partial_2 n_{N-1} & \cdots & 1-x_N \partial_{N-1}n_{N-1} & -n_{N-1} \\
\partial_1 \rho  & \partial_2 \rho  & \cdots & \partial_{N-1}\rho & -  n_N
\end{pmatrix},
\end{align*}
the differential $d\varphi_{x_0}$ of $\varphi$ at $x_0$, is equal to the identity. By the inverse function theorem, there exists $U_{x_0}$, a neighborhood of $x_0$ such that $\varphi:B(x_0,R)^+\to \Omega\cap U_{x_0}$, and $\varphi:B(x_0,R)^-\to (\Omega^c\cap U_{x_0}) \setminus \partial \Omega $ is a diffeomorphism.

\begin{definition}\label{extension}
	For all $x\in U_{x_0}$, set
\begin{align*}
\bar u(x) :=u(\varphi \circ \pi \circ \varphi^{-1}(x)).
\end{align*}
\end{definition}

Up to taking $U_{x_0}$, a neighborhood of $x_0$, smaller, the map $\varphi \circ \pi \circ \varphi^{-1}$ is the identity on $U_{x_0}\cap\overline\Omega$ and its image is contained in $\overline\Omega$. 
\begin{definition}\label{Omega_extended}
We denote by $\xi$ the euclidean metric and we set $ g:=(\varphi\circ\pi\circ\varphi^{-1})^\star\xi$.
For any $x_0\in \partial\Omega$ there exists $R_{x_0}>0$ such that the ball $B(x_0,R_{x_0})\subseteq U_{x_0}$. Thus, since $\Omega$ is compact, there exist $\bar R>0$ such that for all $x_0\in \partial\Omega, B(x_0,\bar R)\subseteq U_{x_0}$. 
	\end{definition}
The metric $g\in C^\infty(\R^N\setminus \partial\Omega)$, in particular $g=\xi$ in $\Omega$. This means that the Christoffel symbols $\Gamma_{ij}^k$ of the metric $ g$ are such that $\Gamma_{ij}^k=0$ in $\Omega$. We claim that $g$ is Lipschitz continuous in $U_{x_0}$. Since $\varphi^{-1}$ is a smooth diffeomorphism on $U_{x_0}$, we know that $ g\in C^{0,1}(U_{x_0})$ if and only if $(\varphi\circ \pi)^\star\xi \in C^{0,1}(B(x_0,R_{x_0}))$. To check this we compute the matrix associated with the metric $(\varphi\circ \pi)^\star\xi(x)$ for $x\in \{x \in B(x_0,R_{x_0})\mid x_N\neq0\}$. When $x_N>0$, the Jacobian matrix of $\varphi\circ \pi$ is given by 
\begin{align*}
&J_{\varphi\circ \pi}(x)=J_{\varphi}(x)\\
&=\begin{pmatrix}
1-x_N\partial_1   n_1 & -x_N\partial_2  n_1 & \cdots & -x_N\partial_{N-1}n_1 &-  n_1 \\
-x_N\partial_1  n_2 & 1-x_N\partial_2   n_2 & \cdots & - x_N \partial_{N-1}n_2 &-  n_2 \\
\vdots & \vdots & & \vdots & \vdots  \\
-x_N\partial_1 n_{N-1} & -x_N\partial_2 n_{N-1} & \cdots & 1-x_N \partial_{N-1}n_{N-1} & -n_{N-1} \\
\partial_1 \rho  & \partial_2 \rho  & \cdots & \partial_{N-1}\rho & -  n_N
\end{pmatrix},
\end{align*}
and, when $x_N<0$, by
\begin{align*}
J_{\varphi\circ \pi}(x)&=\\
&\begin{pmatrix}
1+x_N\partial_1   n_1 & x_N\partial_2  n_1 & \cdots & x_N\partial_{N-1}n_1 &  n_1 \\
x_N\partial_1  n_2 & 1+x_N\partial_2   n_2 & \cdots &  x_N \partial_{N-1}n_2 &  n_2 \\
\vdots & \vdots & & \vdots & \vdots  \\
x_N\partial_1 n_{N-1} & x_N\partial_2 n_{N-1} & \cdots & 1+x_N \partial_{N-1}n_{N-1} & n_{N-1} \\
\partial_1 \rho  & \partial_2 \rho  & \cdots & \partial_{N-1}\rho &   n_N
\end{pmatrix}.
\end{align*}

We compute the symmetric matrix associated with the metric $(\varphi\circ \pi)^\star\xi$ in both cases (see Appendix~\ref{Appendix_B} for details). It is smooth outside of $\{x \in B(x_0,R_{x_0})\mid x_N= 0\}$, and taking the limit when $x_N \to 0$, we have
\begin{align}\label{lim_left_right_g}
\nonumber&\lim\limits_{\stackrel{x_N\to 0}{>}} (\varphi\circ \pi)^\star\xi(x) =\lim\limits_{\stackrel{x_N\to 0}{<}} (\varphi\circ \pi)^\star\xi(x)=\\
 &\begin{pmatrix}
	1+(\partial_1\rho)^2 & \partial_1\rho\,\partial_2\rho & \cdots & \partial_1\rho\,\partial_{N-1}\rho & 0 \\ 
	\partial_1\rho\,\partial_2\rho & 1+(\partial_2\rho)^2 & \cdots & \partial_2\rho\,\partial_{N-1} \rho & 0 \\ 
	\vdots & \vdots & & \vdots & \vdots \\
	\partial_1\rho\,\partial_{N-1} \rho & \partial_2\rho\,\partial_{N-1}\rho & \cdots & 1+(\partial_{N-1}\rho)^2 & 0 \\
	0 & 0 & \cdots & 0 & 1+| \nabla \rho |^2
\end{pmatrix}(x',0).
\end{align}

Since $(\varphi\circ \pi)^\star\xi$ is Lipschitz continuous on $B(x_0,R)^+$ and on $B(x_0,R)^-$, and also that it is continuous on $\{x \in B(x_0,R)\mid x_N=0\}$, we conclude that we can extend $(\varphi\circ \pi)^\star\xi$ as a Lipschitz continuous function on $B(x_0,R)$, and the same holds for $g=(\varphi^{-1})^\star(\varphi\circ\pi)^\star\xi=(\varphi\circ\pi\circ\varphi^{-1})^\star\xi$ on $U_{x_0}$. In addition, we get that $\Gamma_{ij}^k \in L^\infty(U_{x_0})$.


In the rest of the article, we use the following convention: when an index variable is repeated as subscript and superscript, we take the sum over all its possible values. For example, operator $\Delta_{g} = \sum\limits_{i,j=1}^N g^{ij}\left(\partial_{ij}-\sum\limits_{k=1}^N\Gamma_{ij}^k\partial_k\right)$ will be noted $\Delta_{g}=g^{ij}(\partial_{ij}-\Gamma_{ij}^k\partial_k) $. 

\begin{lemma}\label{weak_sol}
	Let $x_0 \in \partial \Omega$, $u \in H^1(\Omega)$ be a solution of (\ref{problem}). Let $U_{x_0}$ be a neighborhood of $x_0$. Let $\bar u$, $\bar a$, and $\bar f$ be as in Definition \ref{extension}, then $\bar u$ satisfies $-\Delta_{ g}\bar u+\bar a \bar u = \bar f$, in the weak sense, i.e.
	$$
	\int_{U_{x_0}} \left(\nabla \bar u,\nabla \phi\right)_{ g}dv_{ g}=\int_{U_{x_0}}( \bar f-\bar a\bar u) \phi \,dv_{g} \text{ for all } \phi \in H_0^1(U_{x_0}),
	$$
	where $dv_{ g}$ denotes the Riemannian volume element.
\end{lemma}
\begin{proof}
	Let $\phi \in H_0^1(U_{x_0})$. 
	
	First, since $\bar u(x)=u(x)$ for all $x\in \Omega \cap U_{x_0}$, we have that 
	$$
	g=\xi\text{ and }\int_{\Omega \cap U_{x_0}} (\nabla \bar u,\nabla \phi)_{ g}dv_{ g}=\int_{\Omega \cap U_{x_0}} (\nabla u,\nabla \phi)dx.
	$$
	Now we define $\tilde\pi:B(x_0,R)^- \rightarrow B(x_0,R)^+:(x',x_N) \rightarrow(x',-x_N)$, so that $\tilde\pi^{-1}$ is well-defined, and we set $ \tilde g:=(\varphi\circ\tilde\pi\circ\varphi^{-1})^\star\xi$.
	
	Secondly, on $\Omega^c \cap U_{x_0}$, we have $\tilde g= g$, and
	\begin{align*}
		&\int_{\Omega^c \cap U_{x_0}} \left(\nabla \bar u,\nabla \phi\right)_{g}dv_{g} = \int_{\Omega^c \cap U_{x_0}} \left(\nabla \bar u,\nabla \phi\right)_{\tilde g}dv_{\tilde g} \\
		=& \int_{B(x_0,R)^-} \left(\nabla  (u\circ\varphi\circ\tilde\pi ),\nabla (\phi\circ\varphi)\right)_{\varphi^\star \tilde g}dv_{\varphi^\star \tilde g}\\
		=& \int_{B(x_0,R)^+} (\nabla  (u\circ\varphi ),\nabla (\phi\circ\varphi\circ\tilde\pi^{-1}))_{(\tilde\pi^{-1})^\star\varphi^\star \tilde g}dv_{(\tilde \pi^{-1})^\star\varphi^\star \tilde g}\\
		=& \int_{\Omega \cap U_{x_0}} (\nabla u ,\nabla (\phi\circ\varphi\circ\tilde\pi^{-1}\circ\varphi^{-1}))_{(\varphi^{-1})^\star(\tilde\pi^{-1})^\star\varphi^\star \tilde g}dv_{(\varphi^{-1})^\star(\tilde\pi^{-1})^\star\varphi^\star \tilde g}\\
		=& \int_{\Omega \cap U_{x_0}} (\nabla u ,\nabla (\phi\circ\varphi\circ\tilde\pi^{-1}\circ\varphi^{-1}))dx
	\end{align*}
	because 
	$$
	(\varphi^{-1})^\star(\tilde\pi^{-1})^\star\varphi^\star \tilde g=(\varphi \circ \tilde\pi^{-1} \circ \varphi^{-1})^\star \tilde g=(\varphi \circ \tilde\pi^{-1} \circ \varphi^{-1})^\star(\varphi\circ\tilde\pi\circ\varphi^{-1})^\star\xi=\xi.
	$$
	Notice that $\phi+\phi(\varphi\circ\tilde\pi^{-1}\circ\varphi^{-1}) \in H_0^1(U_{x_0})$. Thus,
	\begin{align*}
		&\int_{U_{x_0}} \left(\nabla \bar u,\nabla \phi\right)_{ g}dv_{ g}=\int_{\Omega \cap U_{x_0}} (\nabla u ,\nabla (\phi+\phi\circ\varphi\circ\tilde\pi^{-1}\circ\varphi^{-1}))dx\\
		=& \int_{\Omega \cap U_{x_0}} (f-au)  (\phi+\phi\circ\varphi\circ\tilde\pi^{-1}\circ\varphi^{-1})dx\\
		=& \int_{\Omega \cap U_{x_0}} (\bar f-\bar a \bar u) \phi \,dv_g + \int_{\Omega^c \cap U_{x_0}} (\bar f-\bar a \bar u) \phi \,dv_g \\
		=&\int_{U_{x_0}}( \bar f-\bar a\bar u) \phi \,dv_{g},
	\end{align*}
	where we used an integration by parts to get the second equality. 
\end{proof}

\begin{lemma}\label{neumann_extension}
Let $u \in C^2(\overline\Omega)$ and $x=(x',x_N) \in \partial \Omega\cap U_{x_0}$. If $\partial_\nu u(x) =0$, then $\partial_\nu \bar u(x) =0$.	
\end{lemma}

\begin{proof}
If $t<0$, we have $\varphi\circ\pi\circ\varphi^{-1}(x+t\nu) = x+t\nu$, thus we consider $t>0$. Also, notice that $\varphi(x',-t)=x+t\nu$.
		
\begin{align*}
	0&=\partial_\nu u(x)=\lim\limits_{\substack{t \to 0 \\ >}}\frac{u(x)-u(x+t\nu)}{t}=\lim\limits_{\substack{t \to 0 \\ >}} \frac{u(\varphi(x',0))-u(\varphi(x',-t))}{t} \\
	&= \lim\limits_{\substack{t \to 0 \\ >}} \frac{u(\varphi\circ\pi(x',0))-u(\varphi\circ\pi(x',t))}{t} \\
	&= \lim\limits_{\substack{t \to 0 \\ >}} \frac{u(\varphi\circ\pi\circ\varphi^{-1}(x))-u(\varphi\circ\pi\circ\varphi^{-1}(x-t\nu))}{t}\\
	&=\lim\limits_{\substack{t \to 0 \\ >}} \frac{\bar u(x)-\bar u(x-t\nu)}{t}=\partial_{-\nu}\bar u(x)=-\partial_\nu \bar u(x).
\end{align*}	
\end{proof}

\begin{lemma}\label{regular_sol_extended}
	Let $x_0 \in \partial \Omega$. Let $p>1$ and $u \in W^{1,p}(\overline\Omega)$ be a solution of (\ref{problem}) with $f,a\in C^{0,\alpha}(\overline \Omega)$ for all $\alpha \in (0,1)$. Let $U_{x_0}$ be a neighborhood of $x_0$ and $\bar u$ as in Definition \ref{extension}, and let $V_{x_0}\subset\subset U_{x_0}$. Then, $\bar u \in C^{2,\alpha}(V_{x_0})$ for all $\alpha \in (0,1)$ and there exists $C>0$ such that, for all $q>1$,
	\begin{align*}
		\Vert u \Vert_{C^{2,\alpha}(V_{x_0}\cap \overline\Omega)}\leq C\left( \Vert u \Vert_{L^{q}(U_{x_0}\cap \Omega)}+\Vert f \Vert_{C^{0,\alpha}(U_{x_0}\cap \overline\Omega)} \right).
	\end{align*} 
\end{lemma}
\begin{proof}
	We use Lemma \ref{weak_sol} on a neighborhood $\tilde U_{x_0}$ of $x_0$ so that we have $\bar u \in H^1(\tilde U_{x_0})\cap L^\infty( \tilde U_{x_0})$ (also due to our choice of $\varphi$, we have $\bar a \in C^{0,\alpha}(\tilde U_{x_0})\cap L^\infty(\tilde U_{x_0})$). 
	Since $ g^{ij}\in C^{0,1}(\tilde U_{x_0})$ and $ \Gamma^k_{ij}, \bar a\in L^{\infty}(\tilde U_{x_0}) $, by standard elliptic regularity theory \cite[Theorem 9.11]{Gilbarg-Trudinger} there exists $U_{x_0}\subset\subset \tilde U_{x_0}$ such that $\bar u \in W^{2,q}(U_{x_0})$ for all $q\in (1,\infty)$. 
	Then, by Sobolev's embeddings, $\bar u \in C^{1,\alpha}(\overline U_{x_0}) $  for all $\alpha \in (0,1)$.
	We use the chart $\varphi$ that straightens the boundary as in definition~\ref{varphi_pi}. Let $0<R'<R$ be such that $\varphi(B(x_0,R))\subseteq \tilde U_{x_0}$, $\varphi(B(x_0,R'))\subseteq U_{x_0}$ and set $ \hat u :B(x_0,R')\rightarrow \R: x \rightarrow \hat u(x)= \bar u \circ\varphi (x)=u\circ \varphi\circ \pi (x)$, $\hat a = \bar a \circ\varphi$, $\hat f = \bar f \circ\varphi$ and $\hat g = \varphi^\star g = (\varphi \circ \pi)^\star\xi $. 
	We write the equation $\Delta_{\hat g}\hat u+\hat a \hat u = \hat f$ as
	$$
	-\hat g^{ij}\partial_{ij}\hat u = \hat f -\hat a \hat u - \hat g^{ij}\hat \Gamma_{ij}^k\partial_k\hat u \text{ in } B(x_0,R'),
	$$
	where $\hat\Gamma_{ij}^k$ are the Christoffel symbols associated with the metric $\hat g$.
	
	Outside the boundary $\{x_N=0\}$, the function $\hat g^{ij}\hat \Gamma_{ij}^k\partial_k\hat u$ is $\alpha$-Hölder continuous. Therefore it is $\alpha$-Hölder continuous if and only if it is continous on $\{x_N=0\}$. We know that $\partial_\nu\hat u(x',0)= \partial_N \hat u(x',0) = 0$, and after some computations (see Appendix~\ref{Appendix_B} for more details), we get that for all $x \in B(x_0,R')$,
	\begin{align}\label{calcul_difference}
	&\lim\limits_{\stackrel{ x_N\to 0}{x_N>0}}\hat g^{ij}\hat \Gamma_{ij}^k\partial_k\hat u(x)-\lim\limits_{\stackrel{ x_N\to 0}{x_N<0}}\hat g^{ij}\hat \Gamma_{ij}^k\partial_k\hat u(x)\\
	=& \frac{4}{\Vert n \Vert}H(x',0)\partial_N\hat u(x',0)= 0,\nonumber
	\end{align} 
	where $n$ is defined in \eqref{vec_normal}, and $H(x',0)$ denotes the mean curvature of $\partial\Omega$ at $(x',0)$.
	Thus, $-\hat g^{ij}\partial_{ij}\hat u =\hat f -\hat a \hat u - \hat g^{ij}\hat \Gamma_{ij}^k\partial_k\hat u \in C^{0,\alpha}(B(x_0,R'))$. Since $\hat g^{ij}\in C^{0,1}(B(0,R'))$, we use standard elliptic regularity theory \cite[Theorem 9.19]{Gilbarg-Trudinger} and \cite[Theorem 6.2]{Gilbarg-Trudinger} in order to get the existence of $V_{x_0}\subset \subset U_{x_0}$, and $0<R''<R'$ such that $\varphi(B(x_0,R''))\subseteq V_{x_0}$ and such that $\hat u \in C^{2,\alpha}(B(x_0,R''))$ with 
	\begin{align*}
		\Vert \hat u \Vert_{C^{2,\alpha}(B(x_0,R''))}\leq C\left( \Vert \hat u \Vert_{L^{q}(B(x_0,R'))}+\Vert \hat f \Vert_{C^{0,\alpha}(B(x_0,R'))} \right)
	\end{align*}
	for some $C>0$. Getting back to $u$, we have that $ u \in C^{2,\alpha}(V'_{x_0}\cap \overline\Omega)$ and
	\begin{align*}
		\Vert u \Vert_{C^{2,\alpha}(V_{x_0}\cap \overline\Omega)}\leq C\left( \Vert u \Vert_{L^{q}(U_{x_0}\cap \Omega)}+\Vert f \Vert_{C^{0,\alpha}(U_{x_0}\cap \overline\Omega)} \right).
	\end{align*}
\end{proof}

We finalise the construction of the Green's function in the next proposition.
\begin{proposition}\label{exist}
	There exists a Green's function $G$ associated with problem (\ref{problem}). Moreover, there exists $C>0$, such that 
	\begin{align*} 
	| G(x,y) |+|x-y|| \nabla_y G(x,y) |\leq \frac{C}{|x-y|^{N-2}}
	\end{align*}
	for all $x,y \in  \Omega$, $x\neq y$. Also,  $G_x\in C^{1,\alpha}_{loc}\left(\overline \Omega\setminus \{x\}\right)$ for all $\alpha\in (0,1)$ and $x \in \Omega$.
\end{proposition}

\begin{proof}	
	By Theorem \ref{thm_exist}, for all $x\in \Omega$, there exists $w_x \in H^1(\Omega)$ a weak solution of 
	$$
	\begin{cases}
		-\Delta w_x+a w_x &= 0  \text{ in } \Omega \\
		\partial_\nu w_x &=-\partial_\nu \hat G_x   \text{ on }  \partial \Omega.
	\end{cases}
	$$
	By elliptic regularity theory \cite[Theorem 2.3.3.6]{Grisvard}, $w_x\in C^{1,\alpha}_{loc}(\overline\Omega\setminus\{x\}$ for all $\alpha \in (0,1)$. Also, if we take $K\subset \Omega$ compact and $x\in K$, there exists $C(K)>0$ independent of $x$ such that $\Vert w_x \Vert_{C^{1,\alpha}(\overline\Omega)}\leq C(K)$. 
	
	We set 
	$$ 
	G(x,y) := \hat G(x,y) +w_x(y),
	$$ 
	for all $x,y \in \Omega, x\neq y$, where $\hat G$ is as in Proposition~\ref{construct}. Then, we get that $G_x\in C^{1,\alpha}_{loc}\left(\overline \Omega\setminus \{x\}\right)$ for all $\alpha\in (0,1)$ and $x \in \Omega$. Also, for all $x\in K$, we have 
	\begin{align}\label{estimate_G_K}
		|G_x(y)|\leq C(K)|x-y|^{2-N}
	\end{align}
	for all $y\in  \Omega$, $x\neq y$.
	
	Next, we take $r>0$ such that $B(x,3r)\subseteq \Omega$, and set $\tilde G_x(z):= G(x,x+rz)$ and $\tilde a(z)=a(x+rz)$ for all $z\in B(0,3)\setminus\{0\}$. We have that $-\Delta G_x+  a\,G_x=0$ in $\Omega\setminus \{x\}$, and thus, that $-\Delta \tilde G_x+ r^2 \tilde a\,\tilde G_x=0$ in $B(0,3)\setminus\{0\}$. By standard elliptic theory \cite[Theorem 8.32]{Gilbarg-Trudinger}, there exists $C>0$ such that 
\begin{align}\label{estimate_grad_G_K_1}
	&r\Vert \nabla \tilde G_x\Vert_{L^\infty\left(B\left(0,\frac{3}{2}\right)\setminus \overline B\left(0,\frac{2}{3}\right)\right)}\leq C \Vert \tilde G_x\Vert_{L^\infty \left(B\left(0,2\right)\setminus \overline B\left(0,\frac{1}{2}\right)\right)} \leq C r^{2-N}.
\end{align}
Also, we have $\Vert \nabla G_x\Vert_{L^\infty\left(\Omega\setminus \overline B\left(x,r\right)\right)} \leq C\Vert \nabla G_x\Vert_{L^\infty\left(B\left(x,\frac{3}{2} r\right)\setminus \overline B\left(x,\frac{2}{3}r\right)\right)}$ for some $C>0$.
	
	We compute $-\Delta G_x + a G_x$ in the sense of distributions. Let $\phi \in C^2(\overline \Omega)$, we have
	\begin{align*}
		&\int_\Omega (-\Delta \phi + a \phi) G_x \,dy \\
		= & \phi(x) + \int_{\partial \Omega} (\phi \partial_\nu \hat G_x - \partial_\nu \phi \hat G_x) \,d\sigma + \int_\Omega (-\Delta \phi + a \phi) \,w_x \, dy \\
		= & \phi(x) - \int_{\partial \Omega} \partial_\nu \phi G_x \,d\sigma + \int_\Omega (-\Delta w_x + a w_x) \phi \, dy \\
		= & \phi(x) - \int_{\partial \Omega} \partial_\nu \phi G_x \,d\sigma.
	\end{align*}
	
	Now, for $x\in \Omega$, we claim that $G_x\in L^p(\Omega)$ for all $p\in \left[1,\frac{N}{N-2}\right)$ and that there exists $C(p)$ independent of $x$ such that 
	\begin{align}\label{norme_lp_G}
		\Vert G_x \Vert_{L^p(\Omega)}\leq C(p).
	\end{align}
	Indeed, let us define $q=\frac{p}{p-1}$, fix $\psi\in C_c^\infty(\Omega)$ and let $u\in C^2(\overline \Omega)$ be such that 
	\begin{align*}
		\begin{cases}
			-\Delta u + a u &=\psi \text{ in }\Omega\\
			\partial_\nu u &= 0 \text{ on }\partial\Omega.
		\end{cases}
	\end{align*}
	We have $\displaystyle\int_\Omega G_x \psi d y=u(x)$. Thus, using Sobolev's embedding, and elliptic regularity \cite[Theorems 2.3.3.6 and 2.4.2.7]{Grisvard} we get that there exits $C(q)>0$ such that
	\begin{align*}
		\left | \displaystyle\int_\Omega G_x \psi d y \right|\leq\Vert u\Vert_{L^\infty(\Omega)}\leq C(q) \Vert u \Vert_{W^{2,q}(\Omega)}\leq C(q)\Vert \psi- au \Vert_{L^{q}(\Omega)}
	\end{align*}
	for all $\psi\in C_c^\infty(\Omega)$. Then, by duality we have $G_x\in L^p(\Omega)$ for all $p\in \left(1,\frac{N}{N-2}\right)$, and by Hölder inequality, we have that $G_x \in L^1(\Omega)$ so that the claim is proved. Also, $G_x$ is a Green's function. 	
	
	Thanks to \eqref{estimate_G_K} and \eqref{estimate_grad_G_K_1}, it remains to show the estimates when $x\in \Omega \setminus K$. Let $\delta>0$, $x \in \Omega\setminus K$ and $y\in \overline\Omega$ such that $|x-y|\geq \delta$. We know that $G_x$ verifies 
	\begin{align*}
		\begin{cases}
			-\Delta G_x + a G_x &=0 \text{ in }\Omega\setminus \{x\}\\
			\partial_\nu G_x &= 0 \text{ on }\partial\Omega,
		\end{cases}
	\end{align*}
	and that $\Vert G_x\Vert_{L^p(\Omega\setminus \overline B(x,\delta))}\leq C(\delta)$ for all $p>1$. Then by standard elliptic theory \cite[Theorems 2.3.3.6]{Grisvard}, there exists $C(\delta)>0$ such that
	\begin{align}\label{estimate_G_delta}
		&\Vert G_x\Vert_{C^{1,\alpha}(\overline\Omega\setminus B(x,\delta))}\\
		\nonumber &\leq  C(\delta) \Vert G_x \Vert_{W^{2,p}(\Omega\setminus \overline B(x,\delta))}\leq C(\delta)\Vert aG_x \Vert_{L^{p}(\Omega\setminus \overline B(x,\delta))}.
	\end{align}	

	Now we assume that $x \in \Omega\setminus K$ and $y\in \Omega$ such that $|x-y|< \delta$. Let $x_0\in \partial\Omega$, and $V:=V_{x_0}$ a neighborhood of $x_0$. Let $x\in V\cap \Omega$, and $\bar G_x$ be the extension of $G_x$ as in definition \ref{extension}. Let $x^*:=\varphi\circ\tilde \pi^{-1} \circ \varphi^{-1}(x)\in \Omega^c$. We claim that
	\begin{align*}
		(-\Delta_{g}+\bar a)\bar G_x=\delta_x+\delta_{x^*} \text{ in }V
	\end{align*}
	in the sense of distributions.
	Let $\phi\in C_c^\infty(V)$. Separating $V\cap\Omega$ and $V\cap \Omega^c$, and using a change of variable we get that 
	$$
	\int_V \bar G_x (-\Delta_{g}+\bar a )\phi \,dv_{g}=\int_{V\cap\Omega} G_x (-\Delta+a) (\phi + \phi \circ \varphi \circ \tilde\pi^{-1}\circ \varphi^{-1})dy.
	$$
	Using Lemma~\ref{neumann_extension}, and the properties of $G_x$, we get that $\displaystyle\int_V \bar G_x (-\Delta_{g}+\bar a )\phi \,dv_{g}=\phi(x)+\phi(x^*)$ and the claim is proved.
	
	We fix $z\in V$. We claim thet there exists $H_z:V\setminus\{z\}\to\R$ with the following properties:
	\begin{align}\label{H_z}
	\begin{cases}
		(-\Delta_{g}+\bar a)H_z=\delta_z &\text{ in }D'(V),\\
		|H_z(y)|\leq C|z-y|^{2-N} &\text{ for all }y\in V\setminus\{z\},\\
		H_z\in C^1_{loc}(V\setminus\{z\}). 
	\end{cases}	
	\end{align}
	We prove the claim. We define $r(y):=\sqrt{g_{ij}(z)(y-z)^i(y-z)^j}$ for all $y\in V$. We have $r^{2-N}\in C^\infty(V\setminus\{z\})$. We define $f:=(-\Delta_{  g}+\bar a)r^{2-N}$ on $V\setminus\{z\}$. It follows from the properties of $  g$ that $f\in L^\infty_{loc}(V\setminus\{z\})$. Moreover, computing $-g^{ij}\partial_{ij}(r^{2-N})+g^{ij}\Gamma_{ij}^k\partial_k(r^{2-N})$ yield the existence of $C>0$ such that
	\begin{align*}
		|f(y)|\leq C |z-y|^{1-N}\text{ for all }y\in V\setminus\{z\}.
	\end{align*}
	We compute $(-\Delta_{  g}+\bar a)r^{2-N}$ in the sense of distributions, and obtain
	\begin{align*}
		(-\Delta_{  g}+\bar a)r^{2-N}=f+K_z\delta_z \text{ in }D'(V),
	\end{align*}
	where $K_z:=\displaystyle\int_{\partial B(0,1)}(\nabla r(y+z)^{2-N},\nu(y))_{ g(z)}dv_{  g(z)}(y)$.
	
	We define $h$ such that 
	\begin{align*}
		\begin{cases}
			(-\Delta_{  g}+\bar a)h &=f \text{ in }V,\\
			h &=0 \text{ on }\partial V.
		\end{cases}
	\end{align*}
	Using our estimate on $f$ above, the fact that $ g^{ij}\in C^{0,1}(V)$, $ \Gamma^k_{ij}, \bar a\in L^{\infty}(V) $ and standard elliptic theory \cite[Lemma 9.17]{Gilbarg-Trudinger}, we get that $h\in W^{2,p}(V)\cap W_0^{1,p}(V)$ for all $p\in \left(1, \frac{N}{N-1} \right)$ and $h\in C_{loc}^{1,\alpha}(V\setminus\{z\})$. Also, there exists $C(p)>0$ such that 
	\begin{align}\label{norm_W2p_h}
		\Vert h\Vert_{W^{2,p}(V)}\leq C(p) \text{ for all }p\in \left(1,\frac{N}{N-1}\right).
	\end{align}
	Now we let $\alpha\in (N-3,N-2)$, $\epsilon>0$ small and we define
	\begin{align*}
		h_\epsilon(y):=\epsilon^\alpha h(z+\epsilon y) \text{ and }f_\epsilon(y):=\epsilon^{2+\alpha}f(z+\epsilon y)
	\end{align*}
	for all $y \in B(0,2)\setminus \overline B\left(0,\frac{1}{2}\right)$.
	We then have that
	\begin{align}\label{h_eps_equ}
		(-\Delta_{  g_\epsilon}+\bar a_\epsilon)h_\epsilon=f_\epsilon \text{ in }B(0,2)\setminus \overline B\left(0,\frac{1}{2}\right),
	\end{align}
	where $  g_\epsilon=  g(z+\epsilon y)$ and $\bar a_\epsilon:=\epsilon^{2+\alpha}\bar a(z+\epsilon y)$.
	Since $\alpha>N-3$, we have that 
	\begin{align}\label{f_eps_estimate}
		|f_\epsilon(y)|\leq C\epsilon^{\alpha-(N-3)}|y|^{1-N}\leq 2^{N-1}C
	\end{align}
	for all $y\in B(0,2)\setminus \overline B\left(0,\frac{1}{2}\right)$. 
	 We fix $p:=\frac{N}{\alpha+2}\in \left(\frac{N}{N-1} \right)$ and $q:=\frac{N}{\alpha}$. A change of variable, Sobolev's embedding and (\ref{norm_W2p_h}) yield
	\begin{align}\label{h_eps_Lp}
		\Vert h_\epsilon \Vert_{L^q\left(B(0,2)\setminus \overline B\left(0,\frac{1}{2}\right)\right)}\leq C\Vert h \Vert_{L^q(V)}\leq C \Vert h\Vert_{W^{2,p}(V)}\leq C
	\end{align}
	for all $\epsilon>0$ small. It then follows from (\ref{h_eps_equ}), (\ref{f_eps_estimate}), (\ref{h_eps_Lp}) and \cite[Theorem 8.17]{Gilbarg-Trudinger} that there exists $C>0$ such that
	\begin{align*}
		|h_\epsilon(y)|\leq C \text{ for all }y \in \R^N \text{ such that }|y|=1.
	\end{align*}
	Therefore, coming back to $h$, we get that $|h(y)|\leq C|y-z|^{-\alpha}$ for all $|z-y|=\epsilon$. Since $\epsilon$ can be chose arbitrary small and $h$ is bounded outside $z$, we have that for any $\alpha \in (N-3,N-2)$, there exists $C(\alpha)>0$ such that
	\begin{align*}
		|h(y)|\leq C(\alpha)|z-y|^{-\alpha}
	\end{align*}
	for all $y\in V\setminus\{z\}$.
	
	Now we set $H_z:=\frac{1}{K_z} (r^{2-N}-h)$. It follows from the above estimates that $H_z$ satisfies (\ref{H_z}).

	We define $\mu_x:=\bar G_x-H_x-H_{x^*}$. It follows from what precedes that
	\begin{align*}
	\begin{cases}
		(-\Delta_{  g}+\bar a)\mu_x=0 &\text{ in }D'(V),\\
		\Vert \mu_x\Vert_{L^p(V)}\leq C(p) &\text{ for all }p\in\left[1,\frac{N}{N-2} \right).
	\end{cases}
	\end{align*}
	Thus, since $ g^{ij}\in C^{0,1}(V)$ and $ \Gamma^k_{ij}, \bar a\in L^{\infty}(V) $, by standard elliptic theory \cite[Theorem 9.11]{Gilbarg-Trudinger}, $\mu_x\in W_{loc}^{2,p}(V)$ and by Sobolev's embedding and a classical bootstrap argument we have 
	\begin{align}\label{mu_estimate}
		\Vert \mu_x \Vert_{C_{loc}^{0,1}(V')}\leq C(V,V',p),
	\end{align}
	with $V'\subset\subset V$. 
	
	We can now finish the proof of the pointwise estimate on $G$. Indeed, thanks to (\ref{mu_estimate}) and (\ref{H_z}), we know now that there exists $C:=C(V,V',p)>0$ such that
	\begin{align*}
		|\bar G_x(y)| \leq C + C \left (|x-y|^{2-N}+|x^*-y|^{2-N} \right)
	\end{align*}
	for all $x,y\in V'$, $x\neq y$. Taking $V'$ small enough, we have that $|x-y|\leq 2|x^*-y|$ for all $x,y\in V'\cap \Omega$, and thus, we have that
	\begin{align*}
		| G_x(y)| \leq C |x-y|^{2-N}
	\end{align*}
	for all $x,y\in V'\cap \Omega$, $x\neq y$. Combining this with (\ref{estimate_G_K}) and (\ref{estimate_G_delta}) gives
	$$ 
	| G(x,y) |\leq \frac{C}{|x-y|^{N-2}}
	$$ 
	for all $x,y \in  \Omega$, $x\neq y$.
	
	Finally, we still assume that $x \in \Omega\setminus K$, and consider two cases. In the first case, $d(x,\partial \Omega)$ is small enough to take $V_x$ a neighborhood of $x$ such that $V_x \cap \Omega^c$ is non-empty and such that the extension as in Definition~\ref{extension} makes sense on $V_x$. We also take $s>0$ such that $B(x,3s)\subseteq V_x$. Let $\bar G_x$ and $\bar a$ be as in Definition~\ref{extension}. We set $\tilde G_x(z)= \bar G(x,x+sz)$, and $\tilde a(z)= \bar a(x+ sz)$ for all $z\in B(0,3)\setminus\{0\}$. We recall that $g=(\varphi\circ\pi\circ\varphi^{-1})^\star \xi$ and that we have $-\Delta_g \bar G_x+ \bar a\,\bar G_x=0$ in $V_x\setminus \{x\}$, and thus, that $-\Delta_{\tilde g} \tilde G_x+ s^2 \tilde a\,\tilde G_x=0$ in $B(0,3)\setminus\{0\}$, where $\tilde g(z)=g(x+sz)$. We denote by $\tilde\Gamma^k_{ij}$ the Christoffel symbols associated with $\tilde g$. Since $\tilde g\in C^{0,1}(V_x)$ and $ \tilde\Gamma^k_{ij}, \tilde a\in L^{\infty}(V_x) $, by standard elliptic theory \cite[Theorem 8.32]{Gilbarg-Trudinger}, there exists $C>0$ such that 
	\begin{align*}
	s\Vert \nabla \tilde G_x\Vert_{L^\infty\left(B\left(0,\frac{3}{2}\right)\cap \Omega\setminus \overline B\left(0,\frac{2}{3}\right)\right)}\leq C \Vert  \tilde G_x\Vert_{L^\infty \left(B\left(0,2\right)\cap\Omega\setminus \overline B\left(0,\frac{1}{2}\right)\right)} \leq C s^{2-N}.
\end{align*}
Also, we have that $\Vert \nabla G_x\Vert_{L^\infty\left(\Omega\setminus \overline B\left(x,s\right)\right)}\leq C\Vert \nabla G_x\Vert_{L^\infty\left(B\left(x,\frac{3}{2}s\right)\cap \Omega\setminus \overline B\left(x,\frac{2}{3}s\right)\right)}$ for some $C>0$.

In the second case, we take $s>0$ such that $B(x,3s)\subseteq \Omega$, and apply similar arguments as in the first case. We set $\tilde G_x(z)= G(x,x+sz)$, and $\tilde a(z)= a(x+ sz)$ for all $z\in B(0,3)\setminus\{0\}$. We have $-\Delta  G_x+ a\, G_x=0$ in $\Omega\setminus \{x\}$, and thus, that $-\Delta \tilde G_x+ s^2 \tilde a\,\tilde G_x=0$ in $B(0,3)\setminus\{0\}$. By standard elliptic theory \cite[Theorem 8.32]{Gilbarg-Trudinger}, there exists $C>0$ such that 
	\begin{align*}
	s\Vert \nabla \tilde G_x\Vert_{L^\infty\left(B\left(0,\frac{3}{2}\right)\setminus \overline B\left(0,\frac{2}{3}\right)\right)}\leq C \Vert  \tilde G_x\Vert_{L^\infty \left(B\left(0,2\right)\setminus \overline B\left(0,\frac{1}{2}\right)\right)} \leq C s^{2-N}.
\end{align*}
Again, we also have $\Vert \nabla G_x\Vert_{L^\infty\left(\Omega\setminus \overline B\left(x,s\right)\right)}\leq C\Vert \nabla G_x\Vert_{L^\infty\left(B\left(x,\frac{3}{2}s\right)\setminus \overline B\left(x,\frac{2}{3}s\right)\right)}$ for some $C>0$.
These two cases combined with (\ref{estimate_grad_G_K_1}) and (\ref{estimate_G_delta}) ends the proof. 
\end{proof}

\section{Uniqueness, symmetry, extended control and positivity}

In the following, we'll denote by $G$ the Green's function associated with problem (\ref{problem}). We give a few properties that $G$ possesses.

\begin{proposition}\normalfont{(Uniqueness)}\label{unicite}
	Let $x \in \overline\Omega$, if there exists $G_1, G_2 \in L^1(\Omega )$ such that $ \int_\Omega G_i(-\Delta \phi + a \phi ) \, dy = \phi (x)$ for $i\in \{1,2\}$ and all $\phi \in C^2(\overline\Omega)$ such that $\partial_\nu \phi_{|\partial\Omega} =0$, then $G_1= G_2$ almost everywhere.
\end{proposition}

\begin{proof}
	First, we claim that $ G_i \in L^q(\Omega)$ for $i\in\{1,2\}$ and all $q \in \left[1,\frac{N}{N-2}\right)$. Indeed, let $\psi \in C^\infty_c(\Omega)$, and $\phi \in C^2(\overline\Omega)$ such that 
	\begin{align*}
	\begin{cases}
		-\Delta \phi + a \phi &= \psi  \text{ in } \Omega \\
		\partial_\nu \phi &= 0   \text{ on } \partial \Omega.
	\end{cases}
	\end{align*}
	We have 
	\begin{align*}
		\left | \int_\Omega \psi \, G_i \,dy\right |  = |\phi(x)|\leq\Vert\phi\Vert_{L^\infty(\Omega)}.
	\end{align*}
	For $q<\frac{N}{N-2}$, we have $\frac{q}{q-1} >\frac{N}{2}$, and, by Sobolev's embeddings, that $W^{2,\frac{q}{q-1}}(\Omega) \hookrightarrow C^{r,\alpha}(\overline\Omega)$ where $r=0$ or $1$ and $\alpha>0$ are such that $r+\alpha=2-\frac{N}{\frac{q}{q-1}}$. Then by standard elliptic theory \cite[Theorem 2.3.3.6]{Grisvard}, there exists $C>0$ such that
	\begin{align*}
		\Vert\phi\Vert_{L^\infty(\Omega)} \leq C \Vert\phi\Vert_{C^{r,\alpha}(\overline\Omega)}\leq C \Vert\phi\Vert_{W^{2,\frac{q}{q-1}}(\Omega)}\leq C\Vert\psi-a\phi\Vert_{L^{\frac{q}{q-1}}(\Omega)}.
	\end{align*}
	We then get that $\left | \displaystyle\int_\Omega \psi \, G_i \,dy\right | \leq C\Vert\psi\Vert_{L^{\frac{q}{q-1}}(\Omega)}$ for all $\psi \in C^\infty_c(\Omega)$, and $G_i \in L^q(\Omega)$ by duality.
	
	Now, we know that $\displaystyle\int_\Omega \psi (G_1- G_2) \, dy =0$ for all $\psi \in C^\infty_c(\Omega)$, and that $G_1-G_2 \in L^q(\Omega)$ for all $q \in \left[1,\frac{N}{N-2}\right)$, thus, we get that $G_1= G_2$ almost everywhere.
\end{proof}

Next, we extend $G$ up to the boundary and show continuity.

\begin{proposition}\normalfont{(Continuity)} \label{continue}
	We have that $G \in C^0( \overline\Omega \times \overline\Omega \setminus Diag(\overline \Omega))$.
\end{proposition}

\begin{proof}
	We fix $x_\infty\in \partial\Omega$ and $y\in \overline\Omega\setminus\{x_\infty\}$. Let $(x_\alpha)_{\alpha\in \N} \in  \Omega$ such that $\lim\limits_{\alpha\to \infty}x_\alpha =x_\infty$. We claim that we can define 
	$$
	G(x_\infty,y)=\lim\limits_{x_\alpha\to \infty}G(x_\alpha,y)\text{ for all }y\in \overline\Omega\setminus\{x_\infty\},
	$$
	and that $\lim\limits_{\alpha\to\infty}G(x_\alpha,.)=G(x_\infty,.)$ in $C_{loc}^0(\overline\Omega \setminus \{x_\infty \} )$. 
	We set $G(x_\alpha,y) := G_\alpha(y)$, and let $U \subset \subset \overline\Omega\setminus\{x_\infty\} $ an open set. 
	Taking $\alpha$ big enough, we have $x_\alpha \not\in \overline U$. 
	We define the set $$G_U:=\{ G_\alpha | \alpha \in \N \text{ and } x_\alpha \not\in \overline U \} \subset C_{loc}^0(\overline \Omega \setminus\{x_\infty \}).$$ 
	For all $y \in \overline \Omega \setminus\{x_\infty \}$ and for all $G_\alpha \in G_U$, due to Proposition~\ref{exist}, there exists $C>0$ independent of $\alpha$ such that $|G_\alpha(y)|\leq \frac{C}{|x_\alpha -y|^{N-2}}\leq C$. 
	Thus, $G_U$ is uniformly bounded in $\alpha$ on any compact subset of $\overline \Omega \setminus\{x_\infty \}$.
	Notice that $-\Delta G_\alpha = -aG_\alpha$ on $ \Omega \setminus \{x_\alpha \}$. 
	By standard elliptic theory \cite[Lemma 9.11]{Gilbarg-Trudinger}, there exists $C>0$ independent of $\alpha$ such that for all $p>N$,  $$\Vert G_\alpha \Vert_{C^1(\overline U)} \leq C \Vert G_\alpha \Vert_{W^{2,p}( U)} \leq C \Vert G_\alpha \Vert_{L^p( \Omega)} \leq C .$$ 
	Thus $G_U$ is equicontinuous on any compact subset of $\overline \Omega \setminus\{x_\infty \}$. 
	By Arzelà-Ascoli's theorem, there exists $\tilde G$ such that, up to a subsequence, $G_{\alpha} \to \tilde G$ when $\alpha \to \infty$. 
	Since this is true for all $U \subset \subset (\overline\Omega\setminus\{x_\infty\}) $, this convergence happens in $C^0_{loc}(\overline \Omega \setminus \{x_\infty\})$. 
	Also, $\tilde G \in L^1(\Omega)$ thanks to Proposition~\ref{exist} and dominated convergence.
	
	Let $\phi \in C^2(\overline \Omega)$ be such that $\partial_\nu \phi_{|\partial \Omega} = 0$, so that 
	$$\displaystyle\int_\Omega (-\Delta \phi + a \phi ) G_{\alpha} \,dy= \phi(x_{\alpha}). $$ 
	Let $\epsilon >0$, we divide $\Omega$ as  $B(x_\infty, \epsilon)\cup \Omega \setminus \overline B(x_\infty, \epsilon)$, so that
	$$\displaystyle\int_{\overline B(x_\infty, \epsilon)} (-\Delta \phi + a \phi ) G_{\alpha} \,dy + \displaystyle\int_{\Omega \setminus \overline B(x_\infty, \epsilon)} (-\Delta \phi + a \phi ) G_{\alpha} \,dy= \phi(x_{\alpha}) .$$ 
	On one hand, we have 
	$$
	\left | \int_{\overline B(x_\infty, \epsilon)} (-\Delta \phi + a \phi ) G_{\alpha} \,dy  \right | \leq \Vert \phi \Vert_{C^2(\overline \Omega)} \int_{\overline B(x_\infty, \epsilon)} | G_{\alpha}| \,dy,
	$$ 
	and by Proposition~\ref{exist} and dominated convergence,
	$$\lim\limits_{\alpha\to +\infty}G_{\alpha} \in L^1(B(x_\infty, \epsilon)) \text{, so that }\lim\limits_{\epsilon\to 0}\left(\lim\limits_{\alpha\to +\infty}\displaystyle\int_{\overline B(x_\infty, \epsilon)} | G_{\alpha}| \,dy\right) = 0 .
	$$ 
	On the other hand, we have
	\begin{align*}
	&\displaystyle\int_{\Omega \setminus \overline B(x_\infty, \epsilon)} (-\Delta \phi + a \phi ) G_{\alpha} \,dy \underset{\alpha \to +\infty}{\longrightarrow} \displaystyle\int_{\Omega \setminus \overline B(x_\infty, \epsilon)} (-\Delta \phi + a \phi ) \tilde G \,dy, \text{ and }\\ 
	& \displaystyle\int_{\Omega \setminus \overline B(x_\infty, \epsilon)} (-\Delta \phi + a \phi ) \tilde G \,dy\underset{\epsilon \to 0}{\longrightarrow} \displaystyle\int_{\Omega} (-\Delta \phi + a \phi ) \tilde G \,dy.
	\end{align*}
	Thus, by uniqueness of the limit we have$\displaystyle\int_{\Omega} (-\Delta \phi + a \phi ) \tilde G \,dy=\phi(x_\infty)$, and by Proposition~\ref{unicite}, $\tilde G$ does not depend of the choice of the sequence $(x_\alpha)_{\alpha\in \N}$ converging to $x_\infty$. 
	We then set $G_{x_\infty}=\tilde G$, so that the limit $\lim\limits_{x_\alpha\to \infty}G(x_\alpha,y)=G(x_\infty,y)$ makes sense. 
	Furthermore, with $G_{x_\infty}$ now defined, the same arguments imply that $\lim\limits_{\alpha\to\infty}G_{x_\alpha}=G_{x_\infty}$ in $C_{loc}^0(\overline\Omega \setminus \{x_\infty \} )$. Thus the claim is proved. Of course using the same arguments, the claim reamains true when $x_\infty\in \overline\Omega$ and $(x_\alpha)_{\alpha\in \N} \in \overline \Omega$.
	
	Now, let $(x_\alpha)_{\alpha\in \N}, (y_\alpha)_{\alpha\in \N} \in \overline \Omega$ and $x,y \in\overline \Omega, x \neq y$ such that $\lim\limits_{\alpha\to \infty}x_\alpha =x$ and $\lim\limits_{\alpha\to \infty}y_\alpha =y$. Let $\delta >0$ such that $|x-y|> 2 \delta$. Then, since $G_x \in C^1_{loc}(\overline \Omega \setminus \{x\})$, we obtain
		\begin{align*}
			|G(x_\alpha, y_\alpha)-G(x,y)|&\leq  |G_{x_\alpha} (y_\alpha)-G_x(y_\alpha)| + |G_x (y_\alpha)-G_x(y)| \\
			\leq & \Vert G_{x_\alpha}-G_x \Vert_{L^\infty(\overline \Omega \setminus \overline B(x, \delta))} + \Vert G_x \Vert_{C^1(\overline \Omega \setminus \overline B(x, \delta))} |y_\alpha - y|.
		\end{align*}
		We conclude that $\lim\limits_{\alpha \to \infty} G(x_\alpha, y_\alpha)=G(x,y)$.
\end{proof}

A direct consequence is that, now, we can extend the pointwise estimate on $G$. Using the same arguments as in the proof of Proposition~\ref{exist} for $x,y\in \overline \Omega$, we get the following.

\begin{lemma}\label{control_extended}
	There exists $C>0$, such that 
	\begin{align*} 
	| G(x,y) |+|x-y|| \nabla_y G(x,y) |\leq \frac{C}{|x-y|^{N-2}}
	\end{align*}
	for all $x,y \in  \overline\Omega$, $x\neq y$. Also, let $\delta>\delta_0>0$, then $G_x\in C^{1,\alpha}\left(\overline \Omega\setminus B(x,\delta)\right)$ for all $\alpha\in (0,1)$ and $x \in \overline\Omega$.
\end{lemma}

In the case where $a \in C^{0,\alpha}(\overline \Omega)$, we also have the estimate on the second derivative outside of the diagonal.
\begin{lemma}
	Assume $a\in C^{0,\alpha}(\overline \Omega)$ for all $\alpha\in (0,1)$. Let $\delta_0>0$. For all $x,y \in  \overline\Omega$ and $\delta>\delta_0$ such that $|x-y|>\delta$, there exists $C>0$, such that
	\begin{align*}
	|\nabla_y^2 G(x,y) |\leq \frac{C}{|x-y|^{N}}.
	\end{align*}	
\end{lemma}
\begin{proof}
	Let $x \in \overline\Omega$. We know that $-\Delta  G_x+ a\, G_x=0$ in $\Omega\setminus \{x\}$. We set $r=|x-y|$, and define $\tilde \Omega=\{z\in\R^N\mid x+rz\in \Omega\}$, $\tilde G_x(z)= G(x,x+r z)$, and $\tilde a(z)= a(x+ r z)$ for all $z\in \tilde \Omega\setminus\{0\}$. We have $-\Delta \tilde G_x+r^2\tilde a \tilde G_x=0$ in $\tilde\Omega\setminus\{0\}$. By standard elliptic theory \cite[Theorems 6.2 and 9.19]{Gilbarg-Trudinger}, and Lemma~\ref{control_extended}, there exists $C>0$ such that
	
	\begin{align*}
	&r^2\Vert \nabla^2 G_x\Vert_{L^\infty\left(\Omega\setminus \overline B\left(x,\frac{2}{3}r\right)\right)}\leq C \Vert G_x\Vert_{L^\infty \left(\Omega\setminus \overline B\left(x,\frac{r}{2}\right)\right)} \leq C r^{2-N}.
\end{align*} 
\end{proof}

The continuity of $G$ allows us to prove the symmetry property for $G$.

\begin{proposition}\normalfont{(Symmetry)}
	The Green's function is symmetric, i.e.\ $G(x,y)=G(y,x)$ for all $x, y \in \overline\Omega, x\neq y$.
\end{proposition}

\begin{proof}
	Let $f,g\in C_c^\infty(\overline \Omega) $, and let $F(x):= \displaystyle\int_\Omega G(y,x) f(y)\,dy$ for all $x\in \overline \Omega$. Since $G\in L^1(\Omega)$, and $G \in C^0(\overline \Omega \times \overline \Omega \setminus Diag(\overline \Omega))$ by Proposition~\ref{continue}, we have $F \in C^0(\overline \Omega)$.
	
	Let $\phi,\psi \in C^\infty(\overline \Omega)$ be such that
	\begin{align*}
		\left\{
    \begin{array}{ll}
        -\Delta \phi + a \phi  &= f  \text{ in } \Omega \\
        \partial_\nu \phi  &= 0  \text{ on } \partial \Omega
    \end{array}
\right. \qquad \text{and} \qquad \left\{
    \begin{array}{ll}
         -\Delta \psi + a \psi  &= g  \text{ in } \Omega \\
        \partial_\nu \psi  &= 0  \text{ on } \partial \Omega.
    \end{array}
\right.
	\end{align*}
	Using Fubini's theorem, we have
	\begin{align*}
		&\int_\Omega F(x)g(x)\,dx= \int_\Omega F(x)(-\Delta + a) \psi(x)\,dx \\
		=& \int_\Omega \int_\Omega G(y,x) f(y)\,dy (-\Delta + a) \psi(x)\,dx \\
		=& \int_\Omega f(y)\int_\Omega G(y,x) (-\Delta + a) \psi(x)\,dx\, dy \\
		=& \int_\Omega (-\Delta + a) \phi(y) \psi(y) \,dy= \int_\Omega \phi(y) (-\Delta + a) \psi(y) \,dy = \int_\Omega \phi(y) g(y) \,dy.
	\end{align*}
	This means that $\displaystyle\int_\Omega (F(x)-\phi(x))g(x)dx=0$ for all $g\in C^\infty(\overline \Omega)$ and since $F, \phi \in C^0(\overline \Omega)$, we have $F(x)=\phi(x)$. Fixing $x \in \overline \Omega$, by definition of $F$ and of the Green's function, we have $\displaystyle\int_\Omega G(y,x) f(y)\,dy= \displaystyle\int_\Omega G(x,y) f(y)\,dy$ for all $f \in C^\infty(\overline\Omega)$, and thus $G(x,y)=G(y,x)$.
	
\end{proof}

We finally prove the global positivity of $G$. Unlike the Dirichlet case, we see that $G_x > 0$ on $\partial\Omega \setminus \{x\}$ for all $x \in \overline\Omega$.

\begin{proposition}\normalfont{(Positivity)}
	For all $x,y \in \overline\Omega, x\neq y$, we have $ G(x,y)>0$.
\end{proposition}

\begin{proof}
	By Lemma~\ref{control_extended}, we know that $G_x(y) \in C^{1,\alpha}(\overline\Omega\setminus B(x,\delta))$ for $\delta>0$, all $x \in \overline\Omega$ and all $\alpha \in (0,1), $ and we know that, for $x\in \Omega$ fixed, $\lim\limits_{y\to x}G_x(y)=+\infty$, thanks to the proof of Proposition~\ref{exist}. 
	
	If we set $(G_x)\_:=\begin{cases}
		-G_x & \text{ if }G_x<0 \\
		0 & \text{ if }G_x\geq0
	\end{cases}$\,, we have $(G_x)_- \in H^1(\Omega\setminus B(x,\delta))$ and $(G_x)_- = 0$ on $\partial B(x,\delta) $ for a $\delta>0$ small enough and for all $x \in \overline\Omega$. Let $\psi \in C_c^\infty(\Omega \setminus \{x\})$, by definition of the Green's function, we have
	$$
	\int_\Omega G_x(-\Delta \psi + a \psi) dx = 0,
	$$
	and integrating by parts,
	$$
	\int_\Omega (\nabla G_x, \nabla \psi) + a G_x \psi) dx = 0.
	$$
	By density of $C^\infty_c(\Omega \setminus B(x,\delta))$ in $H^1_0(\Omega \setminus B(x,\delta))$, we can take $\psi = (G_x)_-$, and thus $\displaystyle\int_\Omega |\nabla (G_x)_-|^2 + a (G_x)_-^2 dx = 0$. Therefore $(G_x)_-=0$ almost everywhere, and by continuity, $G_x \geq0$ on $\overline\Omega \setminus \{x\} $ for all $x \in \overline\Omega$. 
	We apply the maximum principle on $\overline \Omega \setminus B(x,\delta)$ so that if $G_x=0$ somewhere on $\overline \Omega \setminus B(x,\delta) $, then there exists $y_0 \in \partial \Omega$ such that $G_x(y_0)=0$. And by Höpf's lemma, $\partial_\nu G_x(y_0)<0$, which is absurd. 
	Thus,
	$$
	G_x > 0 \text{ on } \overline\Omega \setminus \{x\} \text{ for all } x \in \overline\Omega.
	$$
\end{proof}

\section{Representation formula}
The next result allows us to differentiate a representation formula.
\begin{lemma}\label{derivation_green_formula}
	Let $x_0\in \overline\Omega$, $r>0$, and $u \in C^2(\overline\Omega)$ be a solution of (\ref{problem}). Then
	\begin{align*}
	u(x)=&\int\limits_{B(x_0, 2r)\cap\Omega} f(y) G (x,y)\,dy\\
	&+\int\limits_{\partial B(x_0, 2r)\cap\Omega}[\partial_\nu u(y) G (x,y)-u(y)\partial_\nu G(x,y)] \,d\sigma(y),
	\end{align*}
	for all $x \in (B(x_0, r)\setminus\{x_0\})\cap\Omega$. And this formula can be differentiated:
	\begin{align*}
\partial_{x_i} u(x)=&\int\limits_{B(x_0, 2r)\cap\Omega} f(y)\partial_{x_i} G (x,y)\,dy\\
&+\int\limits_{\partial B(x_0, 2r)\cap\Omega}[\partial_\nu u(y)\partial_{x_i} G (x,y)-u(y)\partial_{x_i}(\partial_\nu G(x,y))]\,d\sigma(y).
\end{align*}
for $i=1,\cdots,N$.
\end{lemma}

\begin{proof}
	We apply the formula in \cite[Lemma 1.5.3.2]{Grisvard} to $u$ in order to get the identity we would like to differentiate. Let $x \in (B(x_0, r)\setminus\{x_0\})\cap\Omega$. Let $\delta>0$, we write the domain $B(x_0, 2r)\cap\Omega$ as $( B(x_0, 2r)\setminus B(x, \delta)) \cup B(x, \delta)\cap\Omega$. For $y\in (B(x_0, 2r)\setminus B(x, \delta))\cap\Omega$, thanks to Lemma~\ref{control_extended}, we know that there exists $C>0$ such that $\partial_{x_i} G(x,y)\leq\frac{C}{\delta^{N-1}} $. Thus, since $f\in C^0(\overline\Omega)$, we have 
	\begin{align*}
		&\partial_{x_i}\left(\int\limits_{(B(x_0, 2r)\setminus B(x, \delta))\cap\Omega} f(y) G (x,y)\,dy\right)\\
		=&\int\limits_{(B(x_0, 2r)\setminus B(x, \delta))\cap\Omega} f(y) \partial_{x_i} G (x,y)\,dy.
	\end{align*}
	Independently, let $e_i$ a unit vector pointing in the i-th direction, $t\in [0,1]$, and $h>0$ small enough such that $|x+t\,he_i-y|\geq \frac{1}{2}|x-y|$. Using the fundamental theorem of calculus, we have
	\begin{align*}
		&\frac{G(x+he_i,y)-G (x,y)}{h} = \frac{1}{h}\int_0^1 \frac{d}{dt}(G(x+t\,he_i,y))dt\\
		\leq \, & \frac{1}{h}\int_0^1 |h|\left|\nabla G(x+t\,he_i,y)\right|dt\leq \int_0^1 \frac{C}{|x+t\,he_i-y|^{N-1}}dt \\
		\leq\, & \frac{C}{|x-y|^{N-1}}.
	\end{align*}
	
	Now, for $\delta<\frac{|x_0-x|}{2}$, we have that
	\begin{align*}
		&\left|\int\limits_{B(x, \delta)\cap\Omega} f(y)\, \frac{G(x+he_i,y)-G (x,y)}{h}\,dy \right|\\
		\leq &\, C \, \Vert f \Vert_{\infty} \int\limits_{B(x, \delta)\cap\Omega} \frac{1}{|x-y|^{N-1}}dy\leq O(\delta)\xrightarrow[\delta\to0]{}0.
	\end{align*}
	For the boundary term, since $\partial_\nu u$ is bounded on $\partial B(x_0, 2r)\cap\Omega$, and since there exists $C>0$ such that $\partial_{x_i} G(x,y) \leq\frac{C}{r^{N-1}}$ for $y\in \partial B(x_0, 2r)\cap\Omega$ and $x\in  B(x_0, r)\cap\Omega$, we have
	\begin{align*}
		&\partial_{x_i}\left(\int\limits_{\partial B(x_0, 2r)\cap\Omega}\partial_\nu u_\epsilon(y) G (x,y)\,d\sigma(y)\right)\\
		=&\int\limits_{\partial B(x_0, 2r)\cap\Omega}\partial_\nu u_\epsilon(y)\partial_{x_i} G (x,y)\,d\sigma(y).
	\end{align*}
	Similarly,
	\begin{align*}
	&\partial_{x_i}\left(\int\limits_{\partial B(x_0, 2r)\cap\Omega}u_\epsilon(y) \partial_\nu G(x,y)\,d\sigma(y)\right)\\
	=&
		\int\limits_{\partial B(x_0, 2r)\cap\Omega}u_\epsilon(y)\partial_{x_i}(\partial_\nu G(x,y))\,d\sigma(y).
	\end{align*}
\end{proof}

\section{Pointwise estimates}
We prove here pointwise estimates for the Green's function and its derivatives when one of its variables is far from the boundary of $\Omega$.
\begin{proposition}\label{estim_int}
	Let $\Gamma(x,y)=\dfrac{1}{(N-2)\omega_{N-1}|x-y|^{N-2}}$. Let also $x\in K\subset \Omega$ a compact set. Then there exists $C=C(K) >0$, such that
\begin{enumerate}[(i)]
	\item $\left|G(x,y)-\Gamma(x,y)\right|\leq C|x-y|^{3-N}$,
	\item $\left| \nabla_y G(x,y)-\nabla_y \Gamma(x,y)\right|\leq C \begin{cases}
		|\ln|x-y|| &\text{ if }N=3\\
		|x-y|^{3-N} &\text{ if }N>3
	\end{cases}$,
\end{enumerate}
for all $y\in \overline\Omega$, $y\neq x$. 
\end{proposition}

\begin{proof}
	Let $x\in K$. As stated in the proofs of Proposition~\ref{construct} and Proposition~\ref{exist},
	$$
	G_x-\Gamma_x= \sum_{i=1}^k \int_\Omega \Gamma_i(x,z)\Gamma(z,y)\,dz+v_x+w_x,
	$$ 
	with 
	\begin{align*}
		\left|\displaystyle\sum\limits_{i=1}^k\displaystyle\int_\Omega \Gamma_i(x,z)\Gamma(z,y)\,dz\right|\leq C |x-y|^{3-N}
	\end{align*}
	for some $C>0$, and all $y\in \overline\Omega$ $x\neq y$. 
	
	Note that by a similar argument as in the proof of Lemma~\ref{derivation_green_formula}, we can differentiate $G_x-\Gamma_x-v_x-w_x$ under the integral sign. 
	Thus, for $i=1,\cdots,N$, we have $\partial_{y_i}\Gamma(z,y)=\frac{(z-y)_i}{(N-2)\omega_{N-1}|z-y|^{N}}$, and by Giraud's lemma (Proposition~\ref{giraud}), 
	\begin{align*}
		\left|\sum_{i=1}^k\displaystyle\int_\Omega \Gamma_i(x,z)\frac{(z-y)_i}{(N-2)\omega_{N-1}|z-y|^{N}}\,dz\right|\leq C \begin{cases}
		|\ln|x-y|| &\text{ if }N=3\\
		|x-y|^{3-N} &\text{ if }N>3
	\end{cases},
	\end{align*}
	for some $C>0$ and all $y\in \overline\Omega$ $x\neq y$. Moreover, there exists $C>0$ and $C(K)>0$ such that $\Vert v_x \Vert_{C^{1}(\overline \Omega)}\leq C$ and $\Vert w_x \Vert_{C^{1}(K)}\leq C(K)$ for all $x\in K$. This proves (i) and (ii).
\end{proof}

\section{Behavior near the boundary}
\label{section_partial_Omega}
In this section we prove pointwise estimates for the Green's function and its derivatives near the boundary of $\Omega$. We will use the following notation. For any $x=(x_1,\dots,x_{N-1},x_N)\in \R^N_+$, we set $x^*=(x_1,\dots,x_{N-1},-x_N)$.
\begin{proposition}\label{estim_bord}
	Let 
	$$\Gamma_\nu(x,y):=\dfrac{1}{(N-2)\omega_{N-1}}\left(\dfrac{1}{|x-y|^{N-2}}+\dfrac{1}{|x^*-y|^{N-2}} \right)$$ 
	the Green's function of operator $-\Delta$ with Neumann boundary condition on $\partial\R_+^N$. Let $x_0\in \partial \Omega$ and $\varphi$ be a chart as in Definition \ref{smooth_subset} such that $\varphi(0)=x_0$, then there exists $C>0$ and $R>0$ such that
\begin{enumerate}[(i)]
	\item $\left|G(\varphi(x),\varphi(y))-\Gamma_\nu(x,y)\right|\leq C\begin{cases}
		|\ln|x-y|| & \text{ if }N=3 \\
		|x-y|^{3-N} & \text{ if } N>3
	\end{cases}$,
	\item $\left| \partial_{y_i} G(\varphi(x),\varphi(y))-\partial_{y_i} \Gamma_\nu(x,y)\right|\leq C \begin{cases}
		\dfrac{|\ln|x-y||}{|x-y|} & \text{ if }N=3 \\
		|x-y|^{2-N} & \text{ if } N>3
	\end{cases}$
\end{enumerate}
for all $x,y \in B(0,R)\cap\overline{\R^N_+}$, $y\neq x$ and $i= 1, \dots, N$.
\end{proposition}

\begin{proof}
By definition of $\Gamma_\nu$, we have in the sense of distributions that
\begin{align*}
\begin{cases}
	-\Delta (\Gamma_\nu)_x &= \delta_x   \text{ in } \R_+^N, \\
	\partial_\nu (\Gamma_\nu)_x &=0  \text{ on } \partial\R_+^N.
\end{cases}
\end{align*}
Let $\xi$ denote the Euclidean metric, let $R>0$ and take $\varphi:=\varphi_0$ a chart as in Definition \ref{smooth_subset}. 
Let $x\in B(0,R)\cap \R_+^N$, and let $u \in C_c^\infty(\overline \Omega)$. We apply the formula in \cite[Lemma 1.5.3.2]{Grisvard} to $u$, and obtain
	\begin{align*}
		&\int\limits_{\varphi(B(0,R))\cap\Omega }G(\varphi(x),y)(-\Delta+a(y))u(y)dy\\
		=& \,u(\varphi(x))+\int\limits_{\partial\varphi(B(0,R))\cap\Omega }\left[u(y)\partial_\nu G(\varphi(x),y)-G(\varphi(x),y) \partial_\nu u(y) \right] d\sigma(y).
	\end{align*}
	Then, with a change of variables,
	\begin{align*}
		&\int\limits_{B(0,R)\cap \R_+^N}G(\varphi(x),\varphi(y))(-\Delta_{\varphi^\star\xi}+a(\varphi(y)))u(\varphi(y)) dv_{\varphi^\star \xi}(y)\\
		&=u(\varphi(x))+\int\limits_{\partial\varphi(B(0,R))\cap\Omega }\left[u(y)\partial_\nu G(\varphi(x),y)-G(\varphi(x),y) \partial_\nu u(y) \right] d\sigma(y),
	\end{align*}
	where $-\Delta_{\varphi^\star\xi}= -(\varphi^\star\xi)^{ij}\partial^2_{ij}+(\varphi^\star\xi)^{ij} \Gamma_{ij}^k\partial_k $, with $\Gamma_{ij}^k$ being the Christoffel symbols associated with the metric $\varphi^*\xi$, and $dv_{\varphi^\star \xi}$ being the Riemannian volume element. 
Thus, setting $H_{x}(y)=G(\varphi(x),\varphi(y))$ and $\tilde a(y) =a(\varphi(y)) $, we have in the sense of distributions that
\begin{align*}
\begin{cases}
	-\Delta_{\varphi^*\xi}H_x+\tilde a H_x &= \delta_x \text{ in }  B(0,R)\cap \R_+^N, \\
	\partial_\nu H_x &=0 \text{ on } B(0,R)\cap \partial\R_+^N.
\end{cases}
\end{align*}
Since $\varphi$ is a smooth diffeomorphism, the $\Gamma_{ij}^k$'s are smooth. We then compute for $y\in B(0,R)\cap \overline{\R_+^N}$,
\begin{align*}
	&[-\Delta +\tilde a](H_x-(\Gamma_\nu)_x)\\
	=&[(-\Delta +\Delta_{\varphi^*\xi})-\Delta_{\varphi^*\xi} +\tilde a]H_x-[-\Delta +\tilde a](\Gamma_\nu)_x.
\end{align*}
We get 
\begin{align*}
&(-\Delta_{\varphi^*\xi}+\tilde a)H_x-(-\Delta (\Gamma_\nu)_x)=\delta_x-\delta_x=0,\\
&-\tilde a(\Gamma_\nu)_x=O\left(\frac{1}{|x-y|^{N-2}} \right),\\
&(-\Delta +\Delta_{\varphi^*\xi})H_x=((\varphi^*\xi)^{ij}-\xi^{ij})\partial^2_{ij}H_x+(\xi^{ij}-(\varphi^*\xi)^{ij})\Gamma_{ij}^k \partial_k H_x\\
&\qquad\qquad\qquad\quad =O\left(\frac{1}{|x-y|^{N-1}} \right),
\end{align*}
so that $[-\Delta +\tilde a](H_x-(\Gamma_\nu)_x)=O\left(\frac{1}{|x-y|^{N-1}} \right)$.
We apply the formula in \cite[Lemma 1.5.3.2]{Grisvard} on $B(0,R)\cap \R_+^N$, and with Giraud's lemma (Proposition~\ref{giraud}), we have that there exists $C>0$ such that
\begin{align*}
	|H_x(y)-\Gamma_\nu(x,y)|\leq C\int\limits_{B(0,R)\cap \R_+^N}\frac{1}{|y-z|^{N-2}}\frac{1}{|x-z|^{N-1}}dz, \\ \leq C \begin{cases}
		|\ln|x-y|| & \text{ if }N=3 \\
		|x-y|^{3-N} & \text{ if } N>3
	\end{cases}\,,
\end{align*}
which proves (i). 

With the same computations as above, if we set $K(x,y):=H_x(y)-\Gamma_\nu(x,y)$, we get that $-\Delta K_x=f_x$, with $f(x,y)=O\left(\frac{1}{|x-y|^{N-1}}\right)$ for all $x,y \in B(0,R)\cap\overline{\R^N_+}$, $y\neq x$. 
Now, we take $r>0$ such that $B(x,3r)\cap\R^N_+\subset B(0,R)\cap\R^N_+$. We set $\tilde K_x(z):=K(x,x+rz)$ and $\tilde f_x(z):=f(x,x+rz)$ for $z\in B(0,2)\cap\tilde\R^N_+\setminus B(0,\frac{1}{2})$, where $\tilde\R^N_+=\{z\in\R^N\mid x+rz \in\R^N_+\}$. Thus, we have $-\Delta \tilde K_x=r^2 f(x,rz)$. By standard elliptic theory \cite[Theorem 8.32]{Gilbarg-Trudinger} and $(i)$, there exists $C>0$ such that 
\begin{align*}
	&r\Vert \nabla K_x\Vert_{L^\infty\left(B(x,2r)\cap\R^N_+\setminus \overline B\left(x,\frac{1}{2}r\right)\right)}\\
	\leq& C (\Vert K_x\Vert_{L^\infty\left(B(x,3 r)\cap\R^N_+\setminus \overline B\left(x,\frac{1}{3}r\right)\right)}+ r^2\Vert f_x\Vert_{L^\infty\left(B(x,3r)\cap\R^N_+\setminus \overline B\left(x,\frac{1}{3}r\right)\right)} \\
	\leq& C \begin{cases}
		|\ln r| & \text{ if } N=3\\
		r^{3-N} & \text{ if } N>3
	\end{cases}.
\end{align*}
Since $\Vert \nabla K_x\Vert_{L^\infty\left(B(0,R)\cap\R^N_+\setminus \overline B\left(x,2r\right)\right)}\leq \Vert \nabla K_x\Vert_{L^\infty\left(B(x,2r)\cap\R^N_+\setminus \overline B\left(x,\frac{1}{2}r\right)\right)}$, we have shown (ii).
\end{proof}

Finally, we can  conclude that Theorem~\ref{theorem_main} is a direct consequence of Propositions \ref{exist}, \ref{estim_int}, \ref{estim_bord}, and Lemma~\ref{control_extended}.

\appendix
\section{}

We recall two results we used several times in this paper. The first result of this section can be found in \cite[Theorem 8.3]{Salsa}.

\begin{theorem}\label{thm_exist}	
Let $\Omega$ be a bounded Lipschitz domain, $f \in L^2(\Omega)$, $g \in L^{2}(\partial\Omega) $ and $a\in L^\infty(\Omega)$ such that the operator $-\Delta+a$ is coercive. Then, there exists $u \in H^1(\Omega)$ a unique weak solution of the problem

\begin{align*}
	\begin{cases}
		-\Delta u + a u &= f  \text{ in } \Omega \\
		\partial_\nu u &= g   \text{ on } \partial \Omega.
	\end{cases}
\end{align*} 
i.e. there exists a unique $u \in H^1(\Omega)$ such that 
$$
\int_\Omega \left(\nabla u , \nabla \psi \right) \,dx + \int_\Omega a u \,\psi \, dx = \int_\Omega f \, \psi \, dx + \int_{\partial\Omega}g \, \psi \,d\sigma,
$$
for all $\psi \in H^1(\Omega)$.
\end{theorem}

The next result can be found in \cite[Lemma 7.5]{Hebey-bleu}.

\begin{proposition}\label{giraud}
	(Giraud's Lemma)
	Let $\Omega$ be an smooth and bounded open subset of $\R^N$, and $X(x,y), Y(x,y)$ two continuous functions on $\overline\Omega \times \overline\Omega \setminus Diag(\overline \Omega)$ such that there exists $C>0$ such that $|X(x,y)|\leq C|x-y|^{\alpha-N}, |Y(x,y)|\leq C|x-y|^{\beta-N} $ for all $x,y \in \overline\Omega, x\neq y$ with $\alpha, \beta \in (0,N)$. Then, the function $ |Z(x,y)| = \int_\Omega X(x,z) Y(z,y) dz $ is continuous for $x \neq y$ and there exists $C>0$ such that 
	$$
	\begin{cases}
		|Z(x,y)| \leq C|x-y|^{\alpha+\beta-N} & \text{ if } \alpha + \beta < N \\
		|Z(x,y)| \leq C (1+ \log|x-y| ) & \text{ if } \alpha + \beta = N\\
		|Z(x,y)| \leq C & \text{ if } \alpha + \beta > N \\
	\end{cases}
	$$
	for all $x,y \in \overline\Omega, x\neq y$. In the last case $Z(x,y)$ is continuous on $\overline\Omega \times \overline\Omega$.
\end{proposition}
\section{}
\label{Appendix_B}
In this section, we detail the computations involved in (\ref{lim_left_right_g}) and (\ref{calcul_difference}). 

Let $\varphi$ and $\pi$ be as in definition~\ref{varphi_pi}. Consider the metric $\hat g := (\varphi\circ \pi)^\star\xi$. Let us decompose $\hat g $ as follows 
\begin{align*}
	\hat g(x) = \begin{cases}
	\hat g^+(x) \text{ if }x_N>0 \\
	\hat g^-(x) \text{ if }x_N<0
\end{cases}.
\end{align*}
Let us compute the symmetric matrices $\hat g^+$ and $\hat g^-$. We obtain that 
\begin{align}\label{calcul_metric}
	&\begin{cases}
	\hat g^+_{ii}(x)&=(1-x_N \partial_i n_i)^2+\sum\limits_{\stackrel{1\leq j\leq N-1}{j\neq i}}(x_N \partial_{i}n_j)^2+(\partial_i \rho)^2 \\
	&\text{ for }i< N,\\
	\hat g^+_{NN}(x)&=1+|\nabla \rho |^2,\\
	\hat g^+_{ij}(x)&=-2x_N\partial_i n_j+ \sum\limits_{k=1}^{N-1} x_N^2 \partial_i n_k \partial_j n_k +\partial_i \rho \partial_j \rho \\
	&\text{ for }i,j<N \text{ and }i\neq j,\\
	\hat g^+_{iN}(x)&=\sum\limits_{j=1}^{N-1}(x_N n_j\partial_i n_j) \text{ for }i< N,
	\end{cases}\nonumber\\
	&\text{ and that }\\
	\nonumber&\begin{cases}
	\hat g^-_{ii}(x)&=(1+x_N \partial_i n_i)^2+\sum\limits_{\stackrel{1\leq j\leq N-1}{j\neq i}}(x_N \partial_{i}n_j)^2 +(\partial_i \rho)^2 \\
	&\text{ for }i< N,\\
	\hat g^-_{NN}(x)&=1+|\nabla \rho |^2,\\
	\hat g^-_{ij}(x)&= 2x_N\partial_i n_j+ \sum\limits_{k=1}^{N-1} x_N^2 \partial_i n_k \partial_j n_k +\partial_i \rho \partial_j \rho \\&\text{ for }i,j<N \text{ and }i\neq j,\\
	\hat g^-_{iN}(x)&=\sum\limits_{j=1}^{N-1}(x_N n_j \partial_i n_j) \text{ for }i< N.
	\end{cases}
\end{align}
Then, $\lim\limits_{x_N\to 0}g^+(x)=\lim\limits_{x_N\to 0}g^-(x)$, and we have shown (\ref{lim_left_right_g}).

Thanks to \eqref{lim_left_right_g}, we can compute $g^{-1}(x',0)$. We obtain
\begin{align}\label{g^{-1}}
	\begin{cases}
		\vspace{1 mm}\hat g^{ii}(x',0)&=\dfrac{1}{1+|\nabla \rho|^2}\left(1+\sum\limits_{\stackrel{1\leq k\leq N-1}{k\neq i}}(\partial_k \rho)^2\right) \text{ for }i< N,\\
	\vspace{2 mm}\hat g^{NN}(x',0)&=\dfrac{1}{1+|\nabla \rho |^2},\\
	\vspace{2 mm}\hat g^{ij}(x',0)&= \dfrac{- \partial_i \rho \partial_j \rho }{1+|\nabla \rho |^2} \text{ for }i,j<N \text{ and }i\neq j,\\
	\hat g^{iN}(x',0)&=0 \text{ for }i< N.
	\end{cases}	
\end{align}

Now, let us show (\ref{calcul_difference}).
Let $\varphi$, $\pi$ as in Definition~\ref{varphi_pi}, and $R>0$ small enough so that $\varphi(x)$ is well defined on $B(x_0,R)$. Take $\hat g$ as above. Let $\hat\Gamma_{ij}^k$ be the Christoffel symbols associated with the metric $\hat g$. Let $u \in C^2(\overline\Omega)$ be a solution of (\ref{problem}), and $\hat u:= u\circ\varphi \circ \pi$. As in the proof of Lemma~\ref{regular_sol_extended}, $\hat u \in C^{1,\alpha}(B(x_0,R))$. We recall that we aim to show 
\begin{align*}
(\ref{calcul_difference}) \quad	\lim\limits_{\stackrel{ x_N\to 0}{x_N>0}}\hat g^{ij}\hat \Gamma_{ij}^k\partial_k\hat u(x)-\lim\limits_{\stackrel{ x_N\to 0}{x_N<0}}\hat g^{ij}\hat \Gamma_{ij}^k\partial_k\hat u(x)= \frac{4}{\Vert n \Vert}H(x',0)\partial_N\hat u(x',0)
\end{align*}
in $B(x_0,R)$, where $n$ is defined in \eqref{vec_normal}, and $H(x',0)$ denotes the mean curvature of $\partial\Omega$ at $(x',0)$. Notice first that 
\begin{align}\label{lim_partial_k_u}
	\lim\limits_{\stackrel{ x_N\to 0}{x_N>0}}\partial_k \hat u(x)=\lim\limits_{\stackrel{ x_N\to 0}{x_N<0}}\partial_k \hat u(x)
\end{align}
since $\hat u\in C^{1,\alpha}(B(x_0,R))$. By (\ref{lim_left_right_g}), the same holds for $\hat g_{ij}$ and thus also for $\hat g^{ij}$, i.e. 
$$\lim\limits_{ x_N\to 0}\hat g^+_{ij}(x)=\lim\limits_{x_N\to 0}\hat g^-_{ij}(x)
$$ 
and 
$$\lim\limits_{ x_N\to 0}(\hat g^+)^{ij}(x)=\lim\limits_{x_N\to 0}(\hat g^-)^{ij}(x).
$$ 

Thus, the non-zero terms in (\ref{calcul_difference}) must come from the $\hat \Gamma_{ij}^k$'s. We recall that 
$$
\hat \Gamma_{ij}^k=\frac{1}{2}\hat g^{kl}(\partial_j \hat g_{il}+\partial_i \hat g_{lj}-\partial_l \hat g_{ij}).
$$
Furthermore, as in (\ref{lim_left_right_g}), we have that 
\begin{align}\label{lim_g=0}
	\lim\limits_{ x_N\to 0}\hat g_{ij}(x)=\lim\limits_{ x_N\to 0}\hat g^{ij}(x)=0 \text{ when }i=N \text{ or }j=N,\text{ and }i\neq j.
\end{align}  
Now, we write
\begin{align*}
	\hat g^{ij}\hat\Gamma_{ij}^k \partial_k \hat u=\hat g^{ij} \frac{1}{2}\hat g^{kl}(\partial_j \hat g_{il}+\partial_i \hat g_{lj}-\partial_l \hat g_{ij})\partial_k \hat u
\end{align*}
and compute all the terms of \eqref{calcul_difference}. We first get the following result implying that non-zero terms in \eqref{calcul_difference} correspond to the terms with $k=N$ in the sum.
\begin{lemma}
	We have 
	$$
	\lim\limits_{\stackrel{ x_N\to 0}{x_N>0}}\hat g^{ij}\sum_{k=1}^{N-1}\hat \Gamma_{ij}^k\partial_k\hat u(x)-\lim\limits_{\stackrel{ x_N\to 0}{x_N<0}}\hat g^{ij}\sum_{k=1}^{N-1}\hat \Gamma_{ij}^k\partial_k\hat u(x)=0.
	$$
\end{lemma}
\begin{proof}
Thanks to \eqref{lim_g=0}, we have
\begin{align*}
\hat g^{ij}\sum_{k=1}^{N-1}\hat\Gamma_{ij}^k \partial_k \hat u(x)	=&\hat g^{ij} \frac{1}{2}\sum_{k=1}^{N-1}\hat g^{kl}(\partial_j \hat g_{il}+\partial_i \hat g_{lj}-\partial_l \hat g_{ij})\partial_k \hat u(x)\\
	=&\hat g^{ij} \frac{1}{2}\sum_{k,l=1}^{N-1}\hat g^{kl}(\partial_j \hat g_{il}+\partial_i \hat g_{lj}-\partial_l \hat g_{ij})\partial_k \hat u(x)+o(1) \\\text{ as }x_N\to 0.
\end{align*}
Thanks to \eqref{lim_g=0} again, we have 
$$\sum_{j=1}^{N-1}\hat g^{Nj}\sum_{k=1}^{N-1}\hat\Gamma_{Nj}^k \partial_k \hat u(x)=o(1)\text{ as }x_N\to 0,$$
and
$$\sum_{i=1}^{N-1}\hat g^{iN}\sum_{k=1}^{N-1}\hat\Gamma_{iN}^k \partial_k \hat u(x)=o(1)\text{ as }x_N\to 0.$$
We thus have
\begin{align*}
	&\hat g^{ij}\hat\Gamma_{ij}^k \partial_k \hat u(x)\\
	=&\sum_{i,j=1}^{N-1}\hat g^{ij} \frac{1}{2}\sum_{k,l=1}^{N-1}\hat g^{kl}(\partial_j \hat g_{il}+\partial_i \hat g_{lj}-\partial_l \hat g_{ij})\partial_k \hat u(x)\\
	&+\hat g^{NN} \frac{1}{2}\sum_{k,l=1}^{N-1}\hat g^{kl}(\partial_N \hat g_{Nl}+\partial_N \hat g_{lN}-\partial_l \hat g_{NN})\partial_k \hat u(x)+o(1) \text{ as }x_N\to 0.
\end{align*}
Thanks to \eqref{calcul_metric}, we can see that 
\begin{align*}
	&\lim\limits_{\stackrel{ x_N\to 0}{x_N>0}}\sum_{l=1}^{N-1}\hat g^{kl}(\partial_j \hat g_{il}+\partial_i \hat g_{lj}-\partial_l \hat g_{ij})(x)\\
	=&\lim\limits_{\stackrel{ x_N\to 0}{x_N<0}}\sum_{l=1}^{N-1}\hat g^{kl}(\partial_j \hat g_{il}+\partial_i \hat g_{lj}-\partial_l \hat g_{ij})(x),
\end{align*}
for all $i,j,k<N$, that
$$
\lim\limits_{\stackrel{ x_N\to 0}{x_N>0}}\sum_{l=1}^{N-1}\hat g^{kl}(\partial_N \hat g_{Nl}+\partial_N \hat g_{lN})(x)=\lim\limits_{\stackrel{ x_N\to 0}{x_N<0}}\sum_{l=1}^{N-1}\hat g^{kl}(\partial_N\hat g_{Nl}+\partial_N \hat g_{lN})(x),
$$
for all $k<N$, and that 
\begin{align*}
	\lim\limits_{\stackrel{ x_N\to 0}{x_N>0}}\sum_{l=1}^{N-1}\hat g^{kl}(-\partial_l \hat g_{NN})(x)=\lim\limits_{\stackrel{ x_N\to 0}{x_N<0}}\sum_{l=1}^{N-1}\hat g^{kl}(-\partial_l \hat g_{NN})(x),
\end{align*}
for all $k<N$.
This combined with \eqref{lim_partial_k_u} ends the proof
\end{proof}
Thus, we assume now that $k=N$. Thanks to \eqref{lim_g=0}, we have
\begin{align}\label{l=N}
	\hat \Gamma_{ij}^N&=\frac{1}{2}\hat g^{Nl}(\partial_j \hat g_{il}+\partial_i \hat g_{lj}-\partial_l \hat g_{ij})\\
	\nonumber&=\frac{1}{2}\hat g^{NN}(\partial_j \hat g_{iN}+\partial_i \hat g_{Nj}-\partial_N \hat g_{ij})+o(1) \text{ as }x_N\to 0.
\end{align}
We have, using \eqref{lim_g=0} again, that 
$$
\sum_{j=1}^{N-1}\hat g^{Nj}\hat \Gamma_{Nj}^N\partial_N \hat u(x)=o(1) \text{ as }x_N\to 0,
$$
and that  
$$
\sum_{i=1}^{N-1}\hat g^{iN}\hat \Gamma_{iN}^N\partial_N \hat u(x)=o(1) \text{ as }x_N\to 0.
$$
We recall, thanks to \eqref{calcul_metric}, that $\partial_N \hat g_{NN}=0$, and using \eqref{l=N}, we have
\begin{align*}
	&\hat g^{NN}\hat \Gamma_{NN}^N\partial_N \hat u(x)\\
	=&\hat g^{NN}\frac{1}{2}\hat g^{Nl}(\partial_N \hat g_{Nl}+\partial_N \hat g_{lN}-\partial_l \hat g_{NN})\partial_N \hat u(x)\\
	=&\hat g^{NN}\frac{1}{2}\hat g^{NN}(\partial_N \hat g_{NN}+\partial_N \hat g_{NN}-\partial_N \hat g_{NN})\partial_N \hat u(x)+o(1)\\
	=&o(1) \text{ as }x_N\to 0.
\end{align*}
We recall, thanks to \eqref{calcul_metric} again, that $\partial_j \hat g_{iN}+\partial_i \hat g_{Nj}=o(1)$ as $x_N\to 0$ for $i,j<N$, and using \eqref{l=N}, we have
\begin{align*}
	\sum_{i,j=1}^{N-1}\hat g^{ij}\hat \Gamma_{ij}^N\partial_N\hat u(x)&=\sum_{i,j=1}^{N-1}\hat g^{ij}\frac{1}{2}\hat g^{Nl}(\partial_j \hat g_{il}+\partial_i \hat g_{lj}-\partial_l \hat g_{ij})\partial_N\hat u(x)\\
	&=\sum_{i,j=1}^{N-1}\hat g^{ij}\frac{1}{2}\hat g^{NN}(\partial_j \hat g_{iN}+\partial_i \hat g_{Nj}-\partial_N \hat g_{ij})\partial_N\hat u(x)+o(1)\\
	&=\sum_{i,j=1}^{N-1}\hat g^{ij}\frac{1}{2}\hat g^{NN}(-\partial_N \hat g_{ij})\partial_N\hat u(x)+o(1)\\
	&=\begin{cases}
		\sum\limits_{i,j=1}^{N-1}\hat g^{ij}\frac{1}{2}\hat g^{NN}2\partial_i n_j\partial_N\hat u(x)+o(1) &\text{ if }x_N>0\\
		\sum\limits_{i,j=1}^{N-1}\hat g^{ij}\frac{1}{2}\hat g^{NN}(-2\partial_i n_j)\partial_N\hat u(x)+o(1) &\text{ if }x_N<0
	\end{cases},
\end{align*}
as $x_N\to 0$.

Thus, we obtain finally that
\begin{align*}
&\lim\limits_{\stackrel{ x_N\to 0}{x_N>0}}\hat g^{ij}\hat \Gamma_{ij}^k\partial_k\hat u(x)-\lim\limits_{\stackrel{ x_N\to 0}{x_N<0}}\hat g^{ij}\hat \Gamma_{ij}^k\partial_k\hat u(x)\\
	=&\lim\limits_{\stackrel{ x_N\to 0}{x_N>0}}\hat g^{ij}\hat \Gamma_{ij}^N\partial_N\hat u(x)-\lim\limits_{\stackrel{ x_N\to 0}{x_N<0}}\hat g^{ij}\hat \Gamma_{ij}^N\partial_N\hat u(x)\\
	=& 2\hat g^{NN}\left(\sum_{i,j=1}^{N-1}\hat g^{ij}\partial_i n_j\right)\partial_N\hat u(x',0).
\end{align*}
Using \eqref{vec_normal} and \eqref{g^{-1}}, we have
\begin{align*}	
	\hat g^{NN}&=\dfrac{1}{1+|\nabla \rho |^2}=\dfrac{1}{\Vert n \Vert^2},\\
	&\text{ and }\\
	\sum_{i,j=1}^{N-1}\hat g^{ij}\partial_i n_j&=\\
	&\dfrac{\sum\limits_{i=1}^{N-1}\left(1+\sum\limits_{\stackrel{1\leq k\leq N-1}{k\neq i}}(\partial_k \rho)^2\right)\partial^2_{ii}\rho-2 \sum\limits_{\stackrel{1\leq i,j\leq N-1}{i\neq j}} \partial_i \rho\, \partial_j \rho\,\partial^2_{ij}\rho}{1+|\nabla \rho |^2}.
\end{align*}
We conclude by computing $H(x',0)= \frac{1}{2} \text{div}(\nu(x',0))$, with $\nu=\dfrac{n}{\Vert n \Vert}$. We get
\begin{align*}
	2 H(x',0)&=\\
	&\dfrac{\sum\limits_{i=1}^{N-1}\left(1+\sum\limits_{\stackrel{1\leq k\leq N-1}{k\neq i}}(\partial_k \rho)^2\right)\partial^2_{ii}\rho-2 \sum\limits_{\stackrel{1\leq i,j\leq N-1}{i\neq j}} \partial_i \rho\, \partial_j \rho\,\partial^2_{ij}\rho}{(1+|\nabla \rho |^2)^{\frac{3}{2}}}.
\end{align*}
This proves (\ref{calcul_difference}).

\bibliographystyle{amsplain}
\bibliography{biblio}

\end{document}